\documentclass[a4paper,leqno]{article}
\usepackage{amsmath,amssymb,latexsym,theorem,curves,epsfig,graphicx,stmaryrd}
\newcommand{\qed}{\hfill$\blacksquare$}

\newcommand{\Z}{{\mathbb Z}}

\newcommand{\Q}{{\mathbb Q}}

\newcommand{\R}{{\mathbb R}}
\newcommand{\C}{{\mathbb C}}

\theoremstyle{change}
\newtheorem{definition}{{\sc Definition.}}[section]
\newenvironment{defn}{\begin{definition}}{\end{definition}
\smallskip}
\newenvironment{defnl}[1]{\begin{definition} \label{#1}}{\end{definition}
\smallskip}
\newtheorem{Definition}[definition]{{\sc Definition}}

{\theoremstyle{changebreak}
\newtheorem{pdefinition}[definition]{{\sc Definition}}
}

{\theoremstyle{changebreak}
\newtheorem{pDefinition}[definition]{{\sc Definition}}
}

\newenvironment{proof}{\noindent {\it Proof.} \rm \hspace{0.2cm}}{\qed}
\newtheorem{nprop}[definition]{{\sc Proposition.}}
\newenvironment{prop}{\begin{nprop}}{\end{nprop}\smallskip}
\newenvironment{propl}[1]{\begin{nprop} \label{#1}}{\end{nprop}\smallskip}
\newtheorem{Nprop}[definition]{{\sc Proposition}}

{\theoremstyle{changebreak}
\newtheorem{dnprop}[definition]{{\sc Proposition.}}
}
\newenvironment{dprop}{\begin{dnprop}}{\end{dnprop}\smallskip}
\newenvironment{dpropl}[1]{\begin{dnprop} \label{#1}}{\end{dnprop}\smallskip}
{\theoremstyle{changebreak}
\newtheorem{dNprop}[definition]{{\sc Proposition}}
}

\newtheorem{nlemma}[definition]{{\sc Lemma.}}
\newenvironment{lemma}{\begin{nlemma}}{\end{nlemma}\smallskip}
\newenvironment{lemmal}[1]{\begin{nlemma} \label{#1}}{\end{nlemma}\smallskip}
\newtheorem{Nlemma}[definition]{{\sc Lemma}}

{\theoremstyle{changebreak}
\newtheorem{dnlemma}[definition]{{\sc Lemma.}}
}

{\theoremstyle{changebreak}
\newtheorem{dNlemma}[definition]{{\sc Lemma}}
}

\newtheorem{ncoro}[definition]{{\sc Corollary.}}
\newenvironment{coro}{\begin{ncoro}}{\end{ncoro}\smallskip}
\newenvironment{corol}[1]{\begin{ncoro} \label{#1}}{\end{ncoro}\smallskip}
\newtheorem{Ncoro}[definition]{{\sc Corollary}}

{\theoremstyle{changebreak}
\newtheorem{dncoro}[definition]{{\sc Corollary.}}
}

{\theoremstyle{changebreak}
\newtheorem{dNcoro}[definition]{{\sc Corollary}}
}

\newtheorem{ntheorem}[definition]{{\sc Theorem.}}

\newenvironment{theoreml}[1]{\begin{ntheorem} \label{#1}}{\end{ntheorem}\smallskip}
\newtheorem{Ntheorem}[definition]{{\sc Theorem}}

{\theoremstyle{changebreak}
\newtheorem{dntheorem}[definition]{{\sc Theorem.}}
}

{\theoremstyle{changebreak}
\newtheorem{dNtheorem}[definition]{{\sc Theorem}}
}

\newtheorem{remark}[definition]{{\sc Remark.}}
\newenvironment{rem}{\begin{remark} \rm}{\end{remark} \smallskip}

\newtheorem{remarks}[definition]{{\sc Remarks.}}

\newtheorem{notation}[definition]{{\sc Notation.}}

{\theoremstyle{changebreak}
\newtheorem{example}[definition]{{\sc Example.}}
}

\newenvironment{exml}[1]{\begin{example} \label{#1} \rm}{\end{example}\smallskip}
{\theoremstyle{changebreak}
\newtheorem{examples}[definition]{{\sc Examples.}}
}
\newenvironment{exms}{\begin{examples} \rm}{\end{examples}
\smallskip}
\newenvironment{exmsl}[1]{\begin{examples} \label{#1} \rm}{\end{examples}
\smallskip}
\newtheorem{Para}[definition]{}

\newenvironment{paral}[1]{\begin{Para} \label{#1} \rm}{\end{Para} \smallskip}
\begin{document}

\title{Conformal properties of harmonic spinors and lightlike
geodesics in signature $(1,1)$}

\author{Frederik Witt \\ Mathematical Institute \\ University of Oxford\\ witt@maths.ox.ac.uk}

\date{}

\maketitle

\begin{abstract}
\noindent We are studying the harmonic and twistor equation on
Lorentzian surfaces, that is a two dimensional orientable manifold
with a metric of signature $(1,1)$. We will investigate the
properties of the solutions of these equations and try to relate the
conformal invariant dimension of the space of harmonic and twistor
spinors to the natural conformal invariants given by the Lorentzian
metric. We will introduce the notion of semi-conformally flat
surfaces and establish a complete classification of the possible
dimensions for this family.

\medskip

\noindent{\it Subj. Class.}: Differential geometry


\end{abstract}

\section{Introduction}

On every Spin manifold we can canonically construct two first order
differential operators, the so--called Dirac and twistor operator. It
is known that the dimension of their resp. kernels -- the space of
harmonic resp. twistor spinors -- is a conformal invariant, that is
invariant under multiplication of the metric used to define these
operators with a smooth and strictly positive function. A natural
question that arises is to know how these dimensions can be expressed
in terms of conformal invariants given by the (pseudo--)Riemannian
metric. This question is particularly interesting in the case of
Riemannian and Lorentzian surfaces since they carry a natural
conformal structure induced by the isothermal charts.

On Riemannian surfaces, harmonic spinors were studied in \cite{hi74}
and \cite{basc92}. The dimension of the space of harmonic spinors
depends essentially -- unlike the dimension of the space of harmonic
forms -- on the conformal class of the metric and the Spin structure
used to define the Dirac operator. Furthermore, the dimension is
bounded. The purpose of this paper is to study the dimension of the
space of harmonic resp.~twistor spinors in the case of an indefinite,
that is a Lorentzian metric,  and to relate them to lightlike vector fields and lightlike geodesics. We will show that:

\begin{itemize}
\item harmonic and twistor spinors can be described in the same way
\item harmonic and twistor spinors are linked to the global behaviour of the lightlike geodesics and the given Spin structure
\item nowhere vanishing harmonic resp.~twistor spinors cause conformal flatness.
\end{itemize}

Though the Lorentzian theory of surfaces looks similar to its
Riemannian analogue, various new phenomena occur. For instance, we
have uncountably many conformal classes of simply connected Lorentzian
surfaces (see \cite{we96}). On the other hand, the torus is the only compact
surface allowed to carry a Lorentzian metric since the Euler number
has to vanish. However, a Lorentzian torus -- unlike its Riemannian
counterpart -- need not be conformally flat as we have
no regularity properties for solutions of hyperbolic partial differential
equations. Further difficulties are caused by the non--compactness of
the isometry group and the indefiniteness of the scalar product on
the spinor bundle. In order to by--pass these problems, different
approaches are carried out: For Lorentzian tori provided with a
left--invariant metric, we can explicitly compute the kernels by
tools developped in \cite{ba81}. Non--conformally flat examples of
Lorentzian tori are given in section \ref{cldiagmetrics}, where we consider
a particular class of metrics for which the resulting partial
differential equations with respect to the trivial Spin structure can
be explicitly solved. We will generalize these examples by the
observation that harmonic and twistor half spinors might be
interpreted as parallel spinors along one family of lightlike
geodesics. We shall introduce {\em semi-conformally flat} Lorentzian
surfaces which are a particular class of time--orientable,
non--conformally flat Lorentzian surfaces (see \ref{scfdef}) for which a classification of the possible dimensions in dependence on
the Spin structure and the global properties of the lightlike
geodesics is achieved. These surfaces can by characterized by the
existence of a divergence-free lightlike vector
field (cf. \ref{hcfchar}). Furthermore, we will rederive some geometric properties of Lorentzian tori shown in \cite{sa97}.

We now want to state our main result. The lightcones in the tangent space induce
two one-dimensional lightlike distributions which according to section 2 may be
labelled unambiguously by $\mathcal{X}$ and $\mathcal{Y}$ provided the
surface is orientable (which we always tacitly assume). Furthermore, a
lightlike vector field is a section either of the $\mathcal{X}$- or of the
$\mathcal{Y}$-distribution. It makes therefore sense to speak of an
$\mathcal{X}$- or $\mathcal{Y}$-flow if the corresponding vector
fields lie in the $\mathcal{X}$- or $\mathcal{Y}$-distribution, of
$\mathcal{X}$- or $\mathcal{Y}$-geodesics or $\mathcal{X}$- or of
$\mathcal{Y}$-conformal flatness (depending on the divergence-free
lightlike vector field to be $\mathcal{X}$ or $\mathcal{Y}$) etc.. By the
classical Poincar\'e-Bendixson theory for ordinary differential equations on a torus we know
that such lightlike $\mathcal{X}$- and $\mathcal{Y}$-integral curves
or ``lines'' are either closed, asymptotic of a closed curve or
dense. These global properties of the lightlike integral curve give raise to the
notion of ``resonant'' and ``non-resonant'' cylinders and tori. On the other hand
the existence of harmonic spinors imposes some extra conditions on the
holonomy of the principal Spin-bundle: If there is a closed
$\mathcal{X}$- and $\mathcal{Y}$-line, then the lift of this line to
the Spin bundle has to be closed as well. This is what we mean by
``$\mathcal{X}$''- resp. ``$\mathcal{Y}$-triviality''. Let $\delta_{\pm}=dim\:
ker(D^{\pm}:\Gamma(S^{\pm})\to\Gamma(S^{\mp}))$ resp. $\tau_{\pm}=dim\:
ker(P^{\pm}:\Gamma(S^{\pm})\to\Gamma(S^{\mp}))$
denote the dimensions of the spaces of positive/negative harmonic
resp. twistor half-spinors. Then we assert
the following to be true (cf. \ref{scfthm}):

\medskip

\noindent {\sc Theorem.} {\it Let $(M^{1+1},g)$ be a compact
  $\mathcal{X}$--conformally flat Lorentzian surface. Then $\delta_+=\tau_-$ and the only possible dimensions for $\delta_+$ are $0,1$ and $+\infty$. These cases are characterized as follows:
\begin{enumerate}
\item $\delta_+\le 1$ if and only if either \begin{itemize}
        \item there exists a dense $\mathcal{X}$--line in which case we have $\delta_+=0$ for the non--trivial Spin structures, or
        \item $M^{1+1}$ is non-resonant or
        \item there exists no $\mathcal{X}$--trivial resonant cylinder on $M^{1+1}$.
        \end{itemize}
Furthermore, we have $\delta_+=1$ for the trivial Spin structure.
\item $\delta_+=+\infty$ if and only if there exists an $\mathcal{X}$--trivial resonant cylinder on $M^{1+1}$. In this case we have $\delta_+=+\infty$ for every Spin structure.
\end{enumerate}

The same conclusion holds for $\mathcal{Y}$ and $\delta_-$ instead of
$\mathcal{X}$ and $\delta_+$, and an analogous assertion can be stated
for twistor spinors.}

\medskip

The question to what extent this result carries over to general
Lorentzian surfaces remains to be settled.

\begin{center}
{\it Acknowledgements}
\end{center}
I wish to express my gratitude to T.Weinstein (Rutgers University) for
having answered my questions about Lorentzian surfaces, M.S\'{a}nchez
and A.Romero (Universidad de Granada) for having send me their work
which was of great use, the SFB 288 ''Differentialgeometrie und
Quantenphysik'' for its kind invitation in December 1999, and of
course, my former supervisor H.Baum (Humboldt--Universit\"{a}t zu
Berlin) for her mathematical and moral support, and who introduced me
to pseudo-Riemannian Spin geometry.

\section{Lorentzian surfaces}\label{lorsurf}

We will give a brief introduction to the theory of Lorentzian surfaces. For details see \cite{we96}.

A Lorentzian surface $(M^{1+1},g)$ is given by a smooth and orientable
two-dimensional manifold provided with an indefinite metric, that is
$TM^{1+1}$ splits into the direct sum of a {\em timelike} bundle $\xi$
and a {\em spacelike} bundle $\eta$. Furthermore, the lightcone
defined by $g$ is built out of two locally integrable {\em lightlike} (or
{\em isotropic}) distributions. We call these distributions $\mathcal{X}$
and $\mathcal{Y}$ according to the following convention: A vector
$v\in TM^{1+1}$ lies in $\mathcal{X}$ if and only if there exists a
further lightlike vector $w$ such that $(v,w)$ is an oriented basis
and $v+w$ is spacelike. This convention is well-defined, and reversing
the orientation interchanges $\mathcal{X}$ with $\mathcal{Y}$. Therefore, we
can assign to any lightlike object the $\mathcal{X}$-
resp. $\mathcal{Y}$-{\em type} and speak of $\mathcal{X}$-
resp. $\mathcal{Y}$-vector fields, curves, geodesics, etc.. We remark
that lightlike vector fields need not exist globally. In fact, their
global existence is equivalent to the existence of a global orthonormal
basis, so that the orthonormal frame bundle over $M^{1+1}$ is
isomorphic to $M^{1+1}\times SO_+(1,1)$ where $SO_+(1,1)$ denotes the
identity component of the isometry group $O(1,1)$. Equivalently, we may assume the
existence of a non-vanishing timelike vector field. Lorentzian
surfaces which admit such vector fields are said to be {\em
  time-orientable}, since they induce an orientation in a timelike
subbundle. Simply connected Lorentzian surfaces are always
time-orientable. For further reference, we introduce the following
notation: If we fix an othonormal basis $s=(s_1,s_2)$, let
$X_s=s_1+s_2$ and $Y_s=-s_1+s_2$ which are $\mathcal{X}$
resp. $\mathcal{Y}$. By a suitable change of the orthonormal basis,
every $\mathcal{X}$- resp. $\mathcal{Y}$-vector field can be written
in this form.

The local integrability of the lightlike distributions guarantees the
existence of {\em isotropic} resp. {\em isothermal} coordinates
$(x,y)$, so we may locally write $g=\lambda^2dxdy$
resp. $g=\lambda^2(-dx^2+dy^2)$ for a smooth $\lambda\not =0$. In particular,
Lorentzian surfaces are locally conformally flat. Recall that two
metrics $g_1$ and $g_2$ on a manifold $M$ are said to be conformally
equivalent if and only if there is a smooth $\lambda>0$ such that
$g_2=\lambda g_1$. In the case where $g_2$ is flat, we say that $g_1$
is {\em conformally flat}. An atlas consisting of isotropic or
isothermal charts defines -- as in the Riemannian case -- a {\em
  conformal structure} on the surface. These correspond bijectively to
conformal classes of Lorentzian metrics. It should be noted, however,
that the corresponding transition functions have no regularity
properties. The two isotropic distributions $\mathcal{X}$ and
$\mathcal{Y}$ are conformal invariants of the Lorentzian surface
$(M^{1+1},g)$. In fact, they determine the conformal class $[g]$. The
maximal integral curves of the $\mathcal{X}$- and $\mathcal{Y}$-vector
fields which we call $\mathcal{X}$- and $\mathcal{Y}$-{\em lines} or {\em null lines} for short, are further conformal
invariants. The $\mathcal{X}$- and $\mathcal{Y}$-lines through $x$
will be denoted by $l_x$ and $m_x$. As we can locally straighten out
the null lines by choosing isotropic null coordinates, only the global
properties of the null lines encode conformal information. In the case
of a simply connected surface, it can be shown that there are no closed null
lines, that two different null lines intersect in at most one point
and that every null line is properly embedded in $M^{1+1}$.

Since the existence of a Lorentzian metric on a compact surface is
equivalent to $\chi(M)=0$ where $\chi$ denotes the Euler
characteristic, every compact Lorentzian surface is diffeomorphic to a
torus. According to the Poincar\'e--Bendixson theory for ordinary
differential equations on the torus, a null line on a compact
Lorentzian surface is either dense, a closed curve homeomorphic to
$S^1$ which cannot be contracted to a point, or an asymptotic of a
closed null line of the same type.

In \cite{sa97}, the explicit behaviour of the null lines and further properties are discussed for metrics of the form
\begin{equation*}
g_{(x_1,x_2)}=E(x_1)dx_1^2+2F(x_1)dx_1dx_2-G(x_1)dx_2^2.
\end{equation*}
Up to a finite covering, all Lorentzian tori with non--trivial isometry group are of that type. If $G\equiv 0$ resp. $|G|>0$, then $g$ is flat resp. conformally flat. We consider  the family of metrics $\mathcal{G'}$ where $G(0)=0$ and $G$ has only isolated zeros in $p_0=0,\: p_1,\ldots,p_{n-1}\in (0,1),\: p_{n+k}=p_k+1$ for all integer $k$. Then $(\R^{1+1},g)$ is incomplete in the three causal senses for all $g\in\mathcal{G}'$, and so is $(T^{1+1},g)$.

In particular, let us consider the two subfamilies
\begin{eqnarray}
        \mathcal{G}_1 & = & \{g\in \mathcal{G}'\;|\;G_{|(0,1)}>0\} \label{san97g1}\\
        \mathcal{G}_2 & = & \{g\in\mathcal{G}'\;|\;G'(p_i)\not=0,\: F(p_i)F(0)>0,0\le i\le n-1\} \label{san97g2}.
\end{eqnarray}

We have the following proposition (see \cite{sa97}):

\begin{propl}{nullcomp}
Let $g\in \mathcal{G}_1\cup \mathcal{G}_2$.
\begin{enumerate}
\item Let $\eta_0=sgn\: F(0)$. Then
\begin{eqnarray*}
        X_1 & = & G\partial_{x_1}+\left(F+\eta_0\sqrt{EG+F^2}\right)\partial_{x_2} \\
        X_2 & = & \partial_{x_1}+\left(\frac{F-\eta_0\sqrt{EG+F^2})}{G}\right)\partial_{x_2}
\end{eqnarray*}
are two linearly independent isotropic vector fields for $g$ (the choice of $\eta_0$ guaranteeing the existence of the limit in $p_i$).
\item The inextendible null geodesics of $X_2$ are complete. Hence there exists incomplete $X_1$--geodesics since $(T^{1+1},g)$ is lightlike incomplete.
\end{enumerate}
\end{propl}

As we remarked above, a Lorentzian torus need not be conformally
flat. In fact, we have the following characterization of conformal
flatness (see \cite{sa97}):

\begin{propl}{confprop}
Let $(M^{1+1},g)$ be a Lorentzian surface.
\begin{enumerate}
\item If there is a nowhere vanishing time- or spacelike conformal
  vector field, then $(M^{1+1},g)$ is conformally flat. The converse
  is true if in addition $M$ is compact.
\item Every conformally flat compact Lorentzian surface is complete.
\end{enumerate}
\end{propl}

We recall that a vector field $K$ is called {\em conformal} if
$\mathcal{L}_Kg=\sigma g$ for a smooth function $\sigma$ (where
$\mathcal{L}$ denotes the Lie derivative).

\section{Pseudo--Riemannian Spin geometry}\label{prspin}

We will give a brief survey of the relevant Spin geometric features we
use in the fourth section. We focus mainly on the signature $(1,1)$. A general
reference is \cite{ba81}.

Let $\left(\R^{p+q},\langle\cdot,\cdot\rangle_{p,q}\right)$ be the
standard pseudo-euclidean vector space of signature $(p,q)$ where $p$
is the dimension of a maximal timelike subspace. We shall always
assume that the $p$ first vectors of an orthonormal basis are
timelike. For $p+q=2m$ we can identify the associated clifford algebra $\mbox{\sl
  Cliff\/}(\R^{p+q},\langle\cdot,\cdot\rangle_{p,q})$ with
$End(\Delta_{p,q})=\C(2^m)$ obtaining thereby an action $\mu$ of
$Cl_{p,q}$ on $\Delta_{p,q}=\C^{2^m}$. This action will be denoted by
$\cdot$, that is $\mu(x,v)=x\cdot v$. An explicit isomorphism in
signature $(1,1)$ is given by extension of the mapping
\begin{equation} \label{cliffrep} e_1\mapsto\left(\begin{array}{cc}
                                                         0 & i \\
                                                        -i & 0
                                                  \end{array}\right)
\;\mbox{and } e_2\mapsto\left(\begin{array}{cc}
                          0 & i \\
                          i & 0
                   \end{array}\right),
\end{equation}

Next, we define the groups $Spin(p,q)$ and $Spin_+(p,q)$. Let $S_{p,q}=\{x\in\R^{p+q}\;|\;\langle x,x\rangle_{p,q}=1\}$ and $H_{p,q}=\{x\in\R^{p+q}\;|\;\langle x,x\rangle_{p,q}=-1\}$.

\begin{defn}
\begin{eqnarray*}
Spin(p,q)   & = & \{x_1\cdot\ldots x_{2l}\;|\;x_i\in H_{p,q}\cup S_{p,q}\}\\
Spin_+(p,q) & = & \{x_1\cdot\ldots x_{2l}\;|\;x_i\in H_{p,q}\cup S_{p,q}\;\mbox{with an even number of timelike factors}\}
\end{eqnarray*}
\end{defn}

In order to deal simultaneously with the pairs $SO(p,q)/Spin(p,q)$ and
$SO_+(p,q)/Spin_+(p,q)$ we write $G(p,q)$ and $\tilde{G}(p,q)$, where
$\tilde{G}(p,q)=Spin(p,q)$ if $G(p,q)=SO(p,q)$ and
$\tilde{G}(p,q)=Spin_+(p,q)$ if $G(p,q)=SO_+(p,q)$. If we use
(\ref{cliffrep}) to represent $Spin(1,1)$, we get the following:

\begin{lemmal}{spinrep}
$$ Spin(1,1)=\left\{g_a=\left(\begin{array}{cc}
                            a & 0 \\
                            0 & \pm\frac{1}{a}
                        \end{array}\right)|\;a\in\R\backslash\{0\}\right\}$$
\end{lemmal}

The volume element $\omega=e_1\cdot e_2$ of $Cl_{1,1}$ defines --
viewed as an endomorphism of $\Delta_{1,1}$ -- a splitting of
$\Delta_{1,1}$ into the direct sum of $\Delta_{1,1}^+=$ Eigenspace of
$\omega$ for $-1=\langle z_1\rangle$ and $\Delta^-_{1,1}=$ Eigenspace of
$\omega$ for $1=\langle iz_2\rangle$ where $(z_1,z_2)$ denotes the
standard basis of $\C^2$ (the sign convention follows \cite{ba81} and
is motivated by the higher dimensional case). Reversing the
orientation interchanges $\Delta_{1,1}^+$ with $\Delta_{1,1}^-$.

Next, we consider Spin structures of $(M^{p+q},g)$, that is reductions
$(Q,f)$ of the $G(p,q)$-frame bundle $P$ to a $\tilde{G}(p,q)$-bundle
$Q$. We have the following criterion for
the existence of such reductions:

\begin{propl}{obspinstruc}
Let $(M^{p+q},g)$ be a connected pseudo--Riemannian manifold and $TM=\xi^p\oplus\eta^q$ a splitting into a time- and spacelike bundle resp. of maximal rank.
\begin{enumerate}
\item $(M^{p+q},g)$ is Spin if and only if $w_2(TM)=w_1^2(\eta)$ where $w_i\in H^i(M,\Z_2)$ denotes the i-th Stiefel--Whitney--class.
\item If $(M^{p+q},g)$ is time--oriented, then the mapping $Spin(M^{p+q},g)\to \pi_1(P,x),\;{(Q,f)}\mapsto f_*\pi_1(Q,y)$ for $y\in f^{-1}(x)$ is injective. In particular, two Spin structures which are isomorphic as a twofold covering of $P$ are isomorphic as Spin structures.
\item If $Spin(M^{p+q},g)\not= \emptyset$, then $card(Spin(M^{p+q},g))=card(H^1(M^{p+q},\Z_2)).$
\end{enumerate}
\end{propl}

For a proof of (i) see \cite{ka68}, for (ii) and (iii) see \cite{ba81}.

In particular, every time--orientable Lorentzian surface admits a Spin
structure. One is explicitly given by $Q_0=M^{1+1}\times Spin_+(1,1)$
and $f_0(x,a)=(x,\lambda(a))$. This Spin structure will be refered to
as the {\em trivial} one; it is unique up to isomorphism if $M$ is simply connected. On the other hand, if $M$ is time--orientable and compact, then $(M^{1+1},g)$ carries four non--isomorphic Spin structures since $H^1(T^{1+1},\Z_2)=\Z_2\oplus\Z_2$.

If $\tilde{g}=\kappa^2 g$ is conformally equivalent to $g$, we can
canonically  associate a Spin structure $(\tilde{Q},\tilde{f})$ with
every Spin structure $(Q,f)$ on $(M^{p+q},g)$: If $\Phi_\kappa:P_g\to
P_{\tilde{g}}$ is the isomorphism defined by
$s=(s_1,\ldots,s_{p+q})\mapsto
\frac{1}{\kappa}s=(\frac{1}{\kappa}s_1,\ldots,\frac{1}{\kappa}s_{p+q})$,
then the subgroup $\left(\Phi_\kappa\circ f\right)_{\#}\pi_1(Q,q)$ in
$\pi_1(P_{\tilde{g}},\Phi_\kappa(f(q)))$ distinguishes by
\ref{obspinstruc} a Spin structure $(\tilde{Q},\tilde{f})$ which can
be shown to be isomorphic with $(Q,f)$.

Now, we fix a Spin structure $(Q,f)$ over $(M^{p+q},g)$. The associated fiber bundle
$$ S=Q\times _{\tilde{G}(p,q)}\Delta_{p,q} $$
is a complex vector bundle of rank $2^{\left[\frac{p+q}{2}\right]}$ which is
called the {\em spinor bundle} associated with $(Q,f)$. The set of
smooth sections of $S$ is denoted by $\Gamma(S)$: its elements are
called {\em spinor fields} or {\em spinors} for short. A spinor
$\varphi$ may be represented by a $\tilde{G}(p,q)$--equivariant
function $\tilde{\varphi}\in
C^\infty(Q,\Delta_{p,q})^{\tilde{G}(p,q)}$, that is,
$\tilde{\varphi}(qg)=g^{-1}\tilde{\varphi}(q)$ for all $q\in Q$ and
$g\in \tilde{G}(p,q)$. One may also think of a spinor as a collection
of local sections $\tilde{s}:U\to Q$ covering $M$ together with a
family of local trivializations $\varphi_{\tilde{s}}\in
C^\infty(U,\Delta_{p,q})$ verifying
$\varphi_{\tilde{s}g}=g^{-1}\varphi_{\tilde{s}}$. Furthermore, we can
consider the sections of
$S^{\pm}=Q\times_{\tilde{G}(1,1)}\Delta_{1,1}^{\pm}$, where the
fiberwise splitting is induced by the volume element of $\mbox{\sl
  Cliff}(T_xM^{1+1},g_x)$, the Clifford algebra generated by
$(T_xM^{1+1},g_x)$. The corresponding sections in $\Gamma(S^{\pm})$ are
called {\em half spinors}. To emphasize the sign, we also speak of
{\em positive} or {\em negative} (half-)spinors. In the subsequent
chapters, we will also use the following representation of half spinors: If
$\varphi\in\Gamma(S^{\pm})$ is a positive resp. a negative spinor, we can
write $\tilde{\varphi}(g)=\tilde{f}^{\pm}(q)u_{\pm 1}$ for
$\tilde{f}^{\pm}\in C^{\infty}(Q,\C)$. Using the representation of
$Spin(1,1)$ of \ref{spinrep}, the transformation law of $\tilde{f}^+$
is given by $\tilde{f}^+(qg_a)=\frac{1}{a}\tilde{f}^+(q)$ and
$\tilde{f}^-(qg_a)=a\tilde{f}^-(q)$. The same holds if one considers
the complex valued functions $f_{\tilde{s}}$ given by the local
trivializations $\varphi_{\tilde{s}}=f_{\tilde{s}}u_{\pm}$.

As for Riemannian Spin bundles the covariant derivative $\nabla^S:\Gamma(S)\to
\Gamma(T^*M\otimes S)$ is induced by the lift to $Q$ of the
Levi--Civita--connection $Z$ in $P$. Fixing a local orthonormal basis
$s=(s_1,s_2)$, we get
\begin{equation}\label{spinder11}
\nabla^S_V\varphi=[\tilde{s},V(\varphi_{\tilde{s}})-\frac{1}{2}g(\nabla^{LC}_Vs_1,s_2)e_1\cdot e_2\cdot \varphi_{\tilde{s}}].
\end{equation}

We also verify the product rule
$$ \nabla^S_V(W\cdot \varphi)=\left(\nabla^{LC}_VW\right)\cdot \varphi+W\cdot\nabla^S_V\varphi. $$

The main difficulty in the pseudo-Riemannian setup is to define a
suitable scalar product on $S$. This can be done as follows: Assume
$(M^{1+1},g)$ to be time-orientable. Let $TM^{1+1}=\xi^1\oplus\eta^1$ be a
splitting into a (now trivial) time- resp. spacelike vector
bundle. Fix orientations in $\xi$ and $\eta$. $P$ can be reduced to the structure group
$K=SO(1)\times SO(1)$, which is maximal compact in $SO_+(1,1)$. The
reduced bundle is given by $P_\xi=\{(s_1,s_2)\:|\:s_1\mbox{ positively
  oriented in }\xi,\: s_2\mbox{ positively
  oriented in }\eta\}$. Then $\tilde{Q}_\xi=f^{-1}(P_{\xi})$ is the
reduction of $Q$ to $\tilde{K}=(Spin_+(1)\times Spin_+(1))/\Z_2$
which is maximal compact in $Spin_+(1,1)$. We have
$S=\tilde{Q}_\xi\times_{\tilde{K}}\Delta_{1,1}$ and
$TM=P_{\xi}\times_K{\R^{1+1}}$. Let
$(\cdot,\cdot)_{\Delta_{1,1}}$ denote the standard
hermitian product on $\Delta_{1,1}$ which is $\tilde{K}$-invariant,
but not $Spin_+(1,1)$-invariant. We can extend this scalar product to
a fiberwise defined scalar product $(\cdot,\cdot)_{\xi}$ on $S$. Now
let $J_{\xi }:S\rightarrow S,\; J_{\xi }\left( \left
    [ \tilde{q},v\right] \right) =\left[ \tilde{q},e_1\cdot
  v\right]$ for $\tilde{q}\in \tilde{Q}_{\xi }$ and the unit vector
field $e_1\in\xi$. We define
$$ <\varphi ,\psi >_x=\left( J_{\xi }\varphi ,\psi \right)_{\xi x}=(e_1\cdot v,w)_{\Delta_{1,1}}, $$
where $\varphi(x)=[\tilde{q},v]$ and $\psi(x)=[\tilde{q},w]$. This is an indefinite, $Spin_+(1,1)$-invariant scalar product on
$S$. Then the formulae
$$ V(<\varphi ,\psi>)= <\nabla^S_V\varphi ,\psi >+ <\varphi ,\nabla^S_V\psi > $$
and
$$ <X\cdot \varphi ,\psi >=<\varphi ,X\cdot
\psi > $$
hold.

Given $(S,\nabla^S)$ over $(M^{1+1},g)$, we can canonically define two
first order differential operators, namely the
Dirac operator $D$ and the twistor operator $P$. In terms of a local
orthonormal basis $s=(s_1,s_2)$, they are given by
$$ D\varphi=-s_1\cdot\nabla^S_{s_1}\varphi +s_2\cdot\nabla^S_{s_2}\varphi\;\mbox{ and } P\varphi=-s_1\otimes\left(\nabla^S_{s_1}\varphi+\frac{1}{2}s_1\cdot D\varphi\right)+s_2\otimes\left(\nabla^S_{s_2}\varphi+\frac{1}{2}s_2\cdot D\varphi\right). $$

In the fourth section, we will deal with the equations $D\varphi=0$
and $P\psi=0$ called the {\em harmonic} resp. the {\em twistor}\/ equation, the latter being equivalent to
$\nabla^S_V\varphi=-\frac{1}{2}V\cdot D\varphi$ for every $V\in\mathfrak{X}(M^{1+1})$. The solutions are refered to as {\em
  harmonic} resp. {\em twistor} spinors. The vector spaces of
harmonic-- resp. twistor spinors will be denoted by $\mathfrak{H}$ and
$\mathfrak{T}$. Furthermore, we will consider
$\mathfrak{H}_{\pm}=\Gamma(S^\pm)\cap\mathfrak{H}$ and analogously
$\mathfrak{T}_{\pm}=\Gamma(S^\pm)\cap\mathfrak{T}$. The superscript
$^0$ denotes the space of harmonic and twistor spinors resp. their
dimension with respect to the trivial Spin structure. We are mainly
interested in the numbers $\delta_{(\pm)}=dim(\mathfrak{H}_{(\pm)})$
and $\tau_{(\pm)}=dim(\mathfrak{T}_{(\pm)})$, since they have the
following well-known property (see \cite{ba99} and \cite{bfgk91}):

\begin{propl}{confspininv}
Let $\tilde{g}=\lambda g$, and $\tilde{\mathfrak{H}}_{(\pm)}$ resp. $\tilde{\mathfrak{T}}_{(\pm)}$ the space of harmonic resp. twistor (half--) spinors with respect to $\tilde{g}$. Then the maps
\begin{enumerate}
\item $\tilde{\varphi}\in \tilde{\mathfrak{H}}_{(\pm)}\mapsto \lambda^{\frac{1}{4}}\tilde{\varphi}\in \mathfrak{H}_{(\pm)}$
\item $\tilde{\psi}\in\tilde{\mathfrak{T}}_{(\pm)}\mapsto \lambda^{-\frac{1}{4}}\tilde{\varphi}\in \mathfrak{T}_{(\pm)}$
\end{enumerate}
are isomorphisms. In particular, $\delta_{(\pm)}$ and $\tau_{(\pm)}$ are conformal invariants.
\end{propl}

\section{Spinor field equations and lightlike geodesics in signature $(1,1)$}

\subsection{Harmonic and twistor spinors}

\begin{dpropl}{hartwichar}
\begin{enumerate}
\item Let $\varphi $ be in $\Gamma \left( S^{+}\right) $ resp. $\Gamma \left(S^{-}\right) $. Then $\varphi $ is harmonic if and only if
\begin{equation*}
\nabla _{X}^{S}\varphi \equiv 0\;\mbox{resp.}\; \nabla _{Y}^{S}\varphi \equiv 0
\end{equation*}
holds for all $\mathcal{X}$--vector fields $X$ resp. $\mathcal{Y}$--vector fields $Y$.
\item Let $\varphi $ be in $\Gamma \left( S^{+}\right) $ resp. $\Gamma \left(S^{-}\right) $. Then $\varphi $ is twistor if and only if
\begin{equation*}
\nabla _{Y}^{S}\varphi \equiv 0\;\mbox{resp.}\; \nabla _{X}^{S}\varphi \equiv 0
\end{equation*}
holds for all $\mathcal{Y}$--vector fields $Y$ resp. $\mathcal{X}$--vector fields $X$.
\end{enumerate}

\begin{proof}
We prove the assertion only for positive spinors, the remaining cases being showed in the same way.

Let $\varphi\in\Gamma(S^+)$ and let $X$ and $Y$ be a $\mathcal{X}$--
resp. $\mathcal{Y}$--vector field which we write $X=X_s=(s_1+s_2)$ and $Y=Y_s=(-s_1+s_2)$.

(i) Using the local expression of $D$, we see that $\varphi$ is
harmonic if and only if
$s_1\cdot\nabla^S_{s_1}\varphi=s_2\cdot\nabla^S_{s_2}\varphi$ which is
equivalent to
$\nabla^S_{s_1}\varphi=\omega\cdot\nabla^S_{s_2}\varphi=-\nabla^S_{s_2}\varphi$
(with the volume element $\omega=s_1\cdot s_2$). Hence $\varphi$ is harmonic if and only if $\nabla^S_{s_1+s_2}\varphi=0$.

(ii) $\nabla^S_{s_i}\varphi=-\frac{1}{2}s_i\cdot D\varphi$ for
$i=1,2$ is equivalent to $\nabla^S_{s_1}\varphi=-\omega\cdot\nabla^S_{s_2}\varphi$ and $\nabla^S_{s_2}=-\omega\cdot \nabla^S_{s_1}\varphi$, hence to $\nabla^S_{-s_1+s_2}\varphi=0$.
\end{proof}
\end{dpropl}

Since $\nabla^S_V\varphi(x)=[q_x,V^*(\tilde{\varphi})(q_x)]$ for all $q_x\in\pi_Q^{-1}(x)$ and $V\in\mathfrak{X}(M)$ (where $V^*$ denotes the horizontal lift of $V$ to Q), proposition \ref{hartwichar} may be restated as follows:

\begin{paral}{const}{\sc Corollary 1.} {\it
\begin{enumerate}
\item
$\mathfrak{H}_+=\{\tilde{\varphi}\in C^{\infty }\left( Q,\Delta_{1,1}^{+}\right)
^{\tilde{G}\left( 1,1\right) }\mid \tilde{\varphi}\text{ is constant along the horizontal lifts of } \mathcal{X}-\text{curves}\}$
\item
$\mathfrak{H}_{-}=\{\tilde{\varphi}\in C^{\infty }\left
  ( Q,\Delta_{1,1}^{-}\right)^{\tilde{G}\left( 1,1\right) }\mid
\tilde{\varphi}\text{ is constant along the horizontal lifts of } \mathcal{Y}-\text{curves}\}$
\item
$\mathfrak{T}_{+}=\{\tilde{\varphi}\in C^{\infty }\left
  ( Q,\Delta_{1,1}^{+}\right)^{\tilde{G}\left( 1,1\right) }\mid
\tilde{\varphi}\text{ is constant along the horizontal lifts of } \mathcal{Y}-\text{curves}\}$
\item
$\mathfrak{T}_{-}=\{\tilde{\varphi}\in C^{\infty }\left
  ( Q,\Delta_{1,1}^{-}\right)^{\tilde{G}\left( 1,1\right) }\mid
\tilde{\varphi}\text{ is constant along the horizontal lifts of } \mathcal{X}-\text{curves}\}$
\end{enumerate}
}
\end{paral}

A further characterization is given by the formula $\nabla^S_{\gamma '(t)}\varphi=\frac{d}{ds}\mathcal{P}^Q_{\gamma:t+s\to t}\varphi(\alpha(t+s))_{|s=0}$ for any smooth curve $\gamma$, where $\mathcal{P}_{\gamma:t+s\to t}^Q$ denotes the parallel transport of $Q$ along $\gamma$ between the fibers $\pi^{-1}_Q(\gamma(t+s))$ and $\pi_Q^{-1}(\gamma(t))$:

\smallskip

\begin{paral}{parallel}{\sc Corollary 2.} {\it
Let $\varphi \in \Gamma \left( S^{+}\right) $. Then $\varphi $
is a positive harmonic spinor if and only if for any $\mathcal{X}$--curve $\alpha$ joining two points $x$ and $y$ in $M^{1+1}$, we have $\varphi(y)=[\mathcal{P}^Q_{\alpha:x\to y}q,v]$ for $\varphi(x)=[q,v]$. Analogous statements hold for $\mathfrak{H}_-,\mathfrak{T}_+$ and $\mathfrak{T}_-$.
}
\end{paral}

As a first application, we note the following:

\begin{propl}{correspondence}
There is a bijective correspondence between the sets
$\{\varphi\in\mathfrak{H}_+\:|\:\varphi(x)\not=0\mbox{ for all }x\}$
and $\{\varphi\in\mathfrak{T}_-\:|\:\psi(x)\not=0\mbox{ for all }x\}$.

\begin{proof}
If $\varphi\in\Gamma(S^+)$ is given by $\tilde{f}^+u_1$ for
$\tilde{f}^+\in C^\infty(Q,\C)$, we can define a twistor spinor
$\psi_{\tilde{f}^+}\in\Gamma(S^-)$ by
$\tilde{\psi}_{\tilde{f}^+}=\frac{1}{\tilde{f}^+}u_{-1}$, since
$\psi_{\tilde{f}^+}(qg_a)=\frac{1}{\tilde{f}^+(qg_a)}=a\frac{1}{\tilde{f}^+(q)}=a\psi_{\tilde{f}^+}(q)=g^{-1}_a\psi_{\tilde{f}^+}(q)$,
so $\frac{1}{\tilde{f}^+}u_{-1}$ defines indeed a
$\tilde{G}(1,1)$--invariant function. Because of $X^*(\frac{1}{\tilde{f}^+})=0$, $\psi=[q,\frac{1}{\tilde{f}^+(q)}u_{-1}]$ defines a negative twistor spinor.
\end{proof}
\end{propl}

\begin{propl}{denseob}
Let $\varphi \in \mathfrak{H}_{+}$. If $\varphi \left( x\right) =0$, then $\varphi _{\mid l_{x}}\equiv 0$. In particular, we have $\delta_+\le 1$ for any Spin structure if there exists a dense null line on $(M^{1+1},g)$. Analogous statements hold for $\mathfrak{H}_-,\mathfrak{T}_+$ and $\mathfrak{T}_-$.

\begin{proof}
The first assertion is a consequence of the above corollaries. Assume
that there is a $x\in M^{1+1}$ with $l_x$ is dense in $M^{1+1}$. Let
$\varphi_1,\varphi_2\in\mathfrak{H}_+$ with $\varphi_1\not\equiv
0$. Pick $c\in \C$ such that $\varphi_2(x)=c\varphi_1(x)$. Hence $(\varphi_2-c\varphi_1)_{|l_x}\equiv 0$, that is $\varphi_2\equiv c\varphi_1$ for continuity reasons.
\end{proof}
\end{propl}

\subsection{Examples}

We now apply the preceding results to compute some explicit
examples. For the sake of simplicity, we will only deal with positive
harmonic spinors, but all examples extend to the remaining cases in an obvious way.

\subsubsection{Simply connected surfaces}

As observed in section \ref{lorsurf}, two different null lines intersect at most
once and closed null lines cannot exist. Furthermore, the frame bundle
$P$ is trivial since $(M^{1+1},g)$ is time--orientable, and the resulting trivial Spin structure is unique up to isomorphism.

\begin{propl}{1con}
If $M^{1+1}$ is simply connected, then $\delta_+=+\infty$.

\begin{proof}
Let $\beta:[0,1]\to M^{1+1}$ be a $\mathcal{Y}$-curve with lift
$\tilde{\beta}$ to $Q$, and let $f_n:[0,1]\to \C$ be a family of
linearly independent smooth functions whose support is strictly
contained in $[0,1]$. Define a $Spin_+(1,1)$--equivariant function
$\tilde{\varphi}_n:Q\to \Delta_{1,1}^+$ by extending
$\tilde{\varphi}_n(\tilde{\beta}(t))=f_n(t)u_1$ first to $Q_{|\beta}$
by the transitive action of $Spin_+(1,1)$ on the fibers, and secondly
to $M^{1+1}$ by parallel transport on $A=\bigcup_{x\in\beta}l_x$ and $\tilde{\varphi}_n\equiv 0$ on $A^c$.
\end{proof}
\end{propl}

This construction depends crucially on the fact that for simply
connected surfaces, the local and the global behaviour of the null
lines are the same. Therefore we can extend local solutions to global
ones. This observation is the key for the construction of harmonic and
twistor spinors in \ref{musurf}: We will link the global behaviour of the null lines to the spinors; by studying the null lines in the large, we will be able to extend local solutions or to find obstructions for doing so.

\subsubsection{Diagonal metrics on Lorentzian tori}\label{cldiagmetrics}

We consider Lorentzian tori $(T^{1+1},g_\lambda)$ whose metric is given by a {\em diagonal metric}
$$ g_\lambda(x_1,x_2)=-\lambda_1^2(x_1,x_2){dx_1}^2+\lambda_2^2(x_1,x_2){dx_2}^2 $$
for $\lambda_1,\lambda_2\not=0$ in $C^{\infty}(\R^2)^{\Z^2}$.

\medskip

\noindent{\bf Left-invariant metrics}

\noindent First we consider the case where $\lambda_1,\lambda_2$ are
constant. Thus $(T^{1+1},g_\lambda)$ may be seen as a Lie group
provided with a left--invariant metric. Since $\pi_1(T^{1+1})$ has no
two--torsions, we may treat the harmonic and the twistor equation for all four Spin structures simultaneously by tools developed in \cite{ba81} which we will briefly sketch:

The problem is to compare two non--isomorphic Spin structures
$(Q_1,f_1)$ and $(Q_2,f_2)$ and to find conditions for $\Gamma(S_1)$
and $\Gamma(S_2)$ to be isomorphic. Let $(M^{p+q},g)$ be a
pseudo--Riemannian Spin manifold of signature $(p,q)$. Let $R'=\{(q_1,q_2)\in Q_1\times
Q_2\;|\;f_1(q_1)=f_2(q_2)\}$. $\Z_2$ acts naturally on each fiber of
$Q_i$, hence on $R'$. The pair $(R,\mu)$, where $R=R'/\Z_2$ and
$\mu:R\to P,\:[q_1,q_2]\mapsto f_1(q_1)$, is called the {\em
deformation}\/ of $(Q_1,f_1)$ and $(Q_2,f_2)$. If $(Q_1,f_1)$ and
$(Q_2,f_2)$ are isomorphic, then $R$ is ismorphic to
$P\times\Z_2$. $G(p,q)\times\Z_2$ acts on $R$ by
$[q_1,q_2].(A,m)=[q_1a,q_2am]$, where $a\in \lambda^{-1}(A)$. This
action is well defined, therefore providing $R$ with the structure of
a $G(p,q)\times \Z_2$--fiber bundle. Next we define the vector bundle
$E=R/G(p,q)\times_{\Z_2}\R$ over $M$. Its complexification $E^{\C}$ is
given by $E^\C=R/G(p,q)\times_{\Z_2}\C$. Let $\tilde{s}_i:U\to Q_i$ be
two local sections and let
$[\tilde{s}]=[(\tilde{s}_1,\tilde{s}_2)]:U\to R$ where $[\;]$ denotes
the equivalence classes in $R$. Let $e\in E^{\C}$. Then $e$ can be represented in the form $e=[\{\tilde{s}\},z]$ ($\{\;\}$ denoting the equivalence classes in $R/G(p,q)$).

\begin{prop}
The map $\beta:S_1\otimes E^\C\to S_2$ defined by $\beta\left([\tilde{s}_1,v]_x\otimes[\{\tilde{s}\},z]_x\right)=[\tilde{s}_2,zv]_x$ is a vector bundle isomorphism.
\end{prop}

Hence, in the case where $E^\C$ is trivial, the spinor bundles $S_1$
and $S_2$ are isomorphic. For instance, this happens if $(Q_1,f_1)$
and $(Q_2,f_2)$ are isomorphic, for $E^\C$ is then isomorphic to $M^{p+q}\times \Z^2$. Thus equivalent Spin structures induce isomorphic spinor bundles. Furthermore, we yield the following

\begin{coro}
On a surface $M$ each two spinor bundles are isomorphic.

\begin{proof}
For the first Chern class of the complexification $E^{\C}$ of the real
line bundle $E$ holds $2c_1(E^\C)=0$. Since $H^2(M,\Z)=0$ or $\Z$
depending on wether or not $M$ is compact, $E^\C$ must be trivial.
\end{proof}
\end{coro}

Next we want to know how the spinor derivative transforms under this isomorphism:

Let $\nabla^{E^\C}$ be the connection induced by the lift to $R$ of
the Levi--Civita connection of $P$. Then one shows that
$\nabla^{E^\C}$ is flat, hence for
$\eta=[\{\tilde{s}\},z]\in\Gamma(E^{\C}_{|U})$ with $z:U\to \C$, we
have $\nabla^{E^\C}_V\eta=[\{\tilde{s}\},V(\eta)]$. Therefore, the following diagramm commutes for every $V\in \mathfrak{X}(M^{p+q})$:

\begin{center}

\begin{minipage}{10cm}

\unitlength0.5cm

\begin{picture}(16,6)
\put(4.5,0){$\Gamma(S_1\otimes E^{\C})$}
\put(8.3,0.2){\vector(1,0){2.6}}
\put(9.3,0.6){$\beta$}
\put(11.2,0){$\Gamma(S_2)$}
\put(12.1,4){\vector(0,-1){3}}
\put(11.2,4.4){$\Gamma(S_2)$}
\put(12.5,2.4){$\nabla_V^{S_2}$}
\put(4.5,4.4){$\Gamma(S_1\otimes E^\C)$}
\put(8.3,4.6){\vector(1,0){2.6}}
\put(9.3,5){$\beta$}
\put(6.4,4.0){\vector(0,-1){3}}
\put(2.5,2.4){$\nabla_V^{S^1}\otimes_{\nabla^{E^\C}_V} 1$}
\put(9.3,2.4){$\circlearrowright$}
\end{picture}

\end{minipage}

\end{center}

\vspace{0.3cm}

where $\nabla^{S_1}_V\otimes_{\nabla^{E^{\C}}_V}1(\varphi\otimes\eta)=\nabla_V^{S_1}\varphi\otimes\eta+\varphi\otimes\nabla_V^{E^\C}\eta$.

If we assume $E^\C$ to be trivial, then we can choose a nowhere vanishing section $e:M^{p+q}\to E^{\C}$. We define the complex-valued form $\omega_e$ by the equation $\nabla^{E^\C}_Ve=\omega_e(V)e$ and $\alpha_e:\Gamma(S_1)\to\Gamma(S_1\otimes E^{\C})$ by $\alpha_e(\varphi)(x)=\varphi_x\otimes e_x$. Then the following diagramm commutes:

\begin{center}

\begin{minipage}{10cm}

\unitlength0.5cm

\begin{picture}(16,6)
\put(8.5,0){$\Gamma(S_1\otimes E^{\C})$}
\put(12.3,0.2){\vector(1,0){2.6}}
\put(13.3,0.6){$\beta$}
\put(15.2,0){$\Gamma(S_2)$}
\put(16.1,4){\vector(0,-1){3}}
\put(15.2,4.4){$\Gamma(S_2)$}
\put(16.5,2.4){$\nabla_V^{S_2}$}
\put(8.5,4.4){$\Gamma(S_1\otimes E^\C)$}
\put(12.3,4.6){\vector(1,0){2.6}}
\put(13.3,5){$\beta$}
\put(10.4,4.0){\vector(0,-1){3}}
\put(0.5,2.4){$\nabla_V^{S^1}+\omega_e(V)$}
\put(13.3,2.4){$\circlearrowright$}
\put(5.7,0.2){\vector(1,0){2.6}}
\put(4.8,4){\vector(0,-1){3}}
\put(5.7,4.6){\vector(1,0){2.6}}
\put(3.8,4.4){$\Gamma(S_1)$}
\put(3.8,0){$\Gamma(S_1)$}
\put(6.5,5){$\alpha_e$}
\put(6.5,0.6){$\alpha_e$}
\put(7.0,2.4){$\nabla^S_V\otimes_{\nabla^{E^{\C}}_V} 1$}
\put(6.4,2.4){$\circlearrowright$}
\end{picture}

\end{minipage}

\end{center}

\vspace{0.5cm}

Let us now consider the special case of a connected Lie group $G$
provided with a left--invariant metric $g$. Let $p:\tilde{G}\to G$ be
the universal cover of $G$. $\pi_1(G)$ acts as a group of deck transformations. Since $g$ is left--invariant, we can trivialize $P$ by choosing $n$ left--invariant vector fields on $G$, i.e. $P= G\times SO_+(p,q)$. Therefore, $Q_0=G\times Spin_+(p,q)$ with $f_0=id\times\lambda$ defines the trivial Spin structure on $G$. The lifts $\tilde{X}_i$ of the vector fields $X_i$ to $\tilde{G}$ are globally left-- and $\pi_1(G)$--invariant vector fields on $\tilde{G}$, hence $\tilde{P}= \tilde{G}\times SO_+(p,q)$ and $\tilde{Q}_0=\tilde{G}\times Spin_+(p,q),\; \tilde{f}_0=id\times\lambda$. We know that $Spin(G,g)\cong Hom(\pi_1(G),\Z_2)$. Let $\chi\in Hom(\pi_1(G),\Z_2)$. $\pi_1(G,g)$ acts on $\tilde{G}\times Spin_+(p,q)$ through $\chi$ by $\omega.(\tilde{g},a)=(\omega.\tilde{g},\chi(\omega)a)$, where $\omega\in\pi_1(G,e)$. Let $Q_\chi=\tilde{G}\times_{[\pi_1(G,e),\chi]}Spin_+(p,q)$ and $f_\chi:Q_\chi\to P,\; [\tilde{g},a]\mapsto [p(\tilde{g}),\lambda(a)]$. Then $(Q_\chi,f_\chi)$ defines a Spin structure and the following propostion holds (cf. \cite{ba81}):

\begin{dprop}
\begin{enumerate}
\item $Spin(G,g)\cong\{(Q_\chi,f_\chi)\;|\;\chi\in Hom(\pi_1(G),\Z_2)\}$
\item The spinor bundle $S_\chi=Q_\chi\times_{Spin_+(p,q)}\Delta_{p,q}$ associated with the Spin structure $(Q_\chi,f_\chi)$ is given by $S_\chi=\tilde{G}\times_\chi\Delta_{p,q}$.
\item The deformation of $(Q_0,f_0)$ and $(Q_\chi,f_\chi)$ is given by $R_\chi=\left(\tilde{G}/ker(\chi)\right)\times SO_+(p,q)$. Furthermore, $E_\chi^{\C}=\left(\tilde{G}/ker(\chi)\right)\times_{\Z_2}\C$.
\end{enumerate}
\end{dprop}

Assume that we have a nowhere vanishing section
$e_\chi\in\Gamma(E^\C)$. Such a section is given by a map
$\epsilon_\chi:\tilde{G}\to \C$ without zeros such that
$\epsilon(\omega.\tilde{g})=\chi(\omega)\epsilon(\tilde{g})$. If for
$\chi\equiv 1$ we have $\epsilon_1\equiv 1$, we can identify $S_1$
with the trivial Spin structure and $\Gamma(S_1)$ with
$C^\infty(G,\Delta_{p,q})$. Let $(g,Id)$ be a global section of
$P\cong G\times SO_+(p,q)$ and let $[\tilde{g},\mathbf{1}]:U\to
Q_\chi$ be a local lift of this section. Then
$\gamma(x)=[\{\tilde{g}(x)\},\mathbf{1}]\in\Gamma(E^\C_{|U})$
corresponds to this section and $e_{\chi}(x)=[\{\tilde{g}(x)\},\epsilon_\chi(\tilde{g}(x))]=\epsilon_\chi(\tilde{g}(x))\gamma(x)$. Since $\nabla^{E^\C}\gamma=0$, we get
\begin{equation*}
\nabla_V^{E^\C_\chi}e_\chi=\left(d\epsilon_\chi\right)(V^*)\gamma = \epsilon^{-1}_\chi\left(d\epsilon_\chi\right)(V^*)e_\chi,
\end{equation*}
hence $\omega_\chi(V)=\epsilon_\chi^{-1}V^*\left(\epsilon_\chi\right), $
where $V^*$ denotes the lift of $V\in \mathfrak{X}(G)$ to $Q_{\chi}$.

Identifying $\Gamma(S_\chi)$ with $\Gamma(S_1)=C^\infty(G,\Delta_{p,q})$ yields
$ \nabla^{S_\chi}_V\varphi=\nabla^{S_1}_V\varphi+\epsilon_\chi^{-1}V^*(\epsilon_\chi)\varphi $
for $\varphi\in C^\infty(G,\Delta_{p,q})$.

For instance, consider the Lorentzian torus
$(T^{1+1},g_\lambda)$. The universal cover is given by $p:\R^2\to
T^{1+1},\; p(x_1,x_2)=\left(e^{2\pi ix_1},e^{2\pi ix_2}\right)$. Then $\pi_1(T)=\Z\oplus\Z$ acts on $\R^2$ by $(z_1,z_2).(x_1,x_2)=(x_1+z_1,x_2+z_2)$. On the other hand, $Hom(\pi_1(T),\Z_2)$ can be identified with $\Z_2\oplus\Z_2=\{(a_1,a_2)\;|\;a_i\in\{\pm1\}\}$ by $\chi_1=\chi(1\oplus 0)=e^{\frac{i\pi}{2}(1-a_1)}$ and $\chi_2=\chi(0\oplus 1)=e^{\frac{i\pi}{2}(1-a_2)}$. We define
$$ \epsilon_a(x_1,x_2)=e^{\frac{i\pi}{2}\left(x_1(1-a_1)+x_2(1-a_2)\right)}. $$
Then $\epsilon_{(1,1)}\equiv 1$ and $\epsilon_\chi\left((z_1,z_2).(x_1,x_2)\right)=\chi_1^{z_1}\chi_2^{z_2}\epsilon_\chi(x_1,x_2)=\chi(z_1,z_2)\epsilon_\chi(x_1,x_2)$. Now $\omega_\chi(X_s)(x)=\epsilon_\chi^{-1}((s_1+s_2)^*)(\epsilon_\chi)(\tilde{x})$, hence
$$ \omega_\chi(X_s)(x)=\frac{i\pi}{2}\left(\frac{1-a_1}{\lambda_1}+\frac{1-a_2}{\lambda_2}\right). $$
Consequently, for a spinor $\tilde{f}u_1\in\Gamma(S_{(a_1,a_2)}^+)\cong C^\infty(T^{1+1},\C u_1)$ to be harmonic, we yield the following equation:
$$ \nabla^{S_1}_X\tilde{f}+\omega(X)\tilde{f}=\frac{1}{\lambda_1}\partial_{x_1}\tilde{f}+\frac{1}{\lambda_2}\partial_{x_2}\tilde{f}+\frac{i\pi}{2}\left(\frac{1-a_1}{\lambda_1}+\frac{1-a_2}{\lambda_2}\right)\tilde{f}=0. $$
Let the developement of $\tilde{f}$ into a Fourier series be given by $\tilde{f}(x_1,x_2)=\sum_{k,l\in\Z}\tilde{f}_{kl}e^{2\pi i(kx_1+lx_2)}$. Then we get the equation
$$ \sum\limits_{k,l\in\Z}\tilde{f}_{kl}\left(\frac{4k-a_1+1}{\lambda_1}+\frac{4l-a_2+1}{\lambda_2}\right)e^{2\pi i(kx_1+lx_2)}=0. $$
Thus $\tilde{f}\in C^\infty(T^{1+1},\C)$ defines a harmonic spinor if and only if
$$ \tilde{f}_{kl}=0\mbox{ or }4k-a_1+1=-\frac{\lambda_1}{\lambda_2}(4l-a_2+1). $$
Since constants are solutions for the trivial Spin structure and
$4k-a_1+1$ and $4l-a_2+1$ are in $\Z$, we finally find:

\begin{center}
\begin{minipage}{8.0cm}
\begin{tabular}{cccc}
$S_{(a_1,a_2)}$ & $\delta_+$ for & $\frac{\lambda_1}{\lambda_2}\in\Q$ & $\frac{\lambda_1}{\lambda_2}\not\in\Q$ \\ \hline
$(+1,+1)$ & & $+\infty$ & $1$ \\
$(+1,-1)$ & & $+\infty$ & $0$ \\
$(-1,+1)$ & & $+\infty$ & $0$ \\
$(-1,-1)$ & & $+\infty$ & $0$ \\
\end{tabular}
\end{minipage}
\end{center}

\medskip

\noindent{\bf Closed metrics}

\noindent Next, let $\lambda_i\in C^\infty(\R^2)^{\Z^2}$ be two periodic
functions satisfying the additional condition:
$$ \partial_{x_2}\lambda_1+\partial_{x_1}\lambda_2=0, $$
that is the form $\lambda=\lambda_1dx_1-\lambda_2dx_2$ is
closed. Therefore we refer to this type of metrics as {\em
  closed}. Fix the orthonormal frame
$s=(\frac{1}{\lambda_1}\partial_{x_1},\frac{1}{\lambda_2}\partial_{x_2})$.
Then $div(X_s)=0$ -- a fact we will reconsider later. In terms of Fourier coefficients, the closedness condition may be restated as
\begin{equation}\label{symlambda}
l\lambda_{1_{kl}}=k\lambda_{2_{kl}},
\end{equation}
where $\lambda_{i_{kl}}$ denotes the kl-th Fourier coefficient of $\lambda_i$. In particular, we have $\lambda_{1_{0l}}=\lambda_{2_{k0}}=0$ for $l$ and $k$ different from $0$. These formulae will prove useful for the subsequent computations. We fix the trivial Spin structure and trivialize with respect to the basis $s$. Using (\ref{spinder11}), we may rephrase the equation $\nabla_{X_s}^S \tilde{f}u_1=0$ as
$$ -\lambda_2\partial_{x_1}\tilde{f}=\lambda_1\partial_{x_2}\tilde{f}+\frac{1}{2}(\partial_{x_2}\lambda_1+\partial_{x_1}\lambda_2), $$
so using our additional assumption gives
\begin{equation}\label{anaeq}
-\lambda_2\partial_{x_1}\tilde{f}=\lambda_1\partial_{x_2}\tilde{f}.
\end{equation}

We remark that the constant spinor $u_1$ defines a solution that
projects onto the torus. Furthermore, if there exists a dense null
line, we already know by \ref{denseob} that there cannot exist any
further linearly independent solutions.

In order to find non-trivial solutions we define the function
$$
\tilde{f}_\alpha(x_1,x_2)=\exp\left(i\pi\alpha(\int_0^{x_1}(\lambda_1(s,x_2)+\lambda_1(s,0))ds-\int_0^{x_2}(\lambda_2(x_1,t)+\lambda_2(0,t))dt)\right)
$$
for $0\not=\alpha\in\R$. Then $\tilde{f}_\alpha$ defines a solution
for (\ref{anaeq}) on $\R^2$. In order to be inducible on the torus
$T$, we must impose the double-periodicity of $\tilde{f}_\alpha$, that
is $\tilde{f}_\alpha(x_1+n,x_2+m)=\tilde{f}_\alpha(x_1,x_2)$ for $n,m\in\Z$. If $l_i=\iint_{T^{1+1}}\lambda_idT^{1+1}$ denotes the 0--th Fourier coefficient of $\lambda_i$, we get the following criterion:

\begin{lemma}
$\tilde{f}_\alpha\in C^\infty(\R^2,\C)^{\Z^2}$ if and only if $\alpha l_i\in \Z$ for $i=1,2$.

\begin{proof}
We have $\tilde{f}_\alpha\left( x_{1}+n,x_{2}+m\right)=\tilde{f}_\alpha\left( x_{1},x_{2}\right)$
if and only if
\begin{equation*}
\alpha \left(n\int\limits_0^1\left(\lambda_1(s,x_2)+\lambda _1(s,0)\right)ds-m\int\limits_0^1\left(\lambda_2(x_1,t)+\lambda _2(0,t)\right)dt\right)\in 2\Z.
\end{equation*}
Now, $\int_{0}^{1}\lambda_{1}(s,x_2)ds=l_{1}$ and
$\int_{0}^{1}\lambda_{2}(x_1,t)ds=l_2$ (use (\ref{symlambda})), so $\tilde{f}_\alpha\in C^\infty(\R^2,\C)^{\Z^2}$ if and only if $\alpha(nl_1-ml_2)\in \Z$ for every $n$ and $m$ in $\Z$.
\end{proof}
\end{lemma}

In particular, if $\frac{l_1}{l_2}=\frac{p}{q}\in\Q$ (where $p$ and $q$ have no common divisor), then $\tilde{f}_{\frac{q}{l_2}}$ defines a solution. Since the set $\{\tilde{f}_{\alpha m}u_1\;|\;m\in\Z\}$ is linearly independent in $\mathfrak{H}_+^0$, we have $\delta_+^0=+\infty$.

\begin{exml}{analex}
Let $c>0$ be a rational number and $f:\R\to\R$ be a smooth function
without zeros with period $1$ such that $-f(2x)\not= c$ for every
$x\in\R$ and $f(2p_i)=-\frac{c}{2}$ for $p_0=0,p_1,\ldots,p_n=1\in [0,1]$, but $f'(2p_i)\not= 0$.
For instance, we could choose $c=2$ and
$f(x)=\frac{1}{10}\cos(2\pi(x+\frac{1}{4}))-1$. Then the diagonal
metric defined by $\lambda_1(x_1,x_2)=-f(x_1-x_2)$ and
$\lambda_2=-f(x_1-x_2)-c$ is closed. Furthermore,
$\frac{l_1}{l_2}=\frac{1}{1-c}\in \Q$, hence $\delta_+^0=+\infty$. If
we express $g_\lambda$ in the new coordinates $(x,y)$ given by
$x=\frac{x_1-x_2}{2}$ and $y=\frac{x_1+x_2}{2}$, we get
$g_\lambda(x,y)=(2cf(2x)+c^2)dx^2-4(f^2(2x)+2f(2x)+2)dxdy+(2cf(2x)+c^2)dy^2$.
Because of our assumptions $g_\lambda\in \mathcal{G}_2$ (see
\ref{nullcomp}), so $g_\lambda$ provides an example of a
non--conformally flat diagonal metric since it is not complete.
\end{exml}

\begin{lemma}
There are no dense $\mathcal{X}$--null lines if and only if $\frac{l_1}{l_2}\in\Q$.

\begin{proof}
As the global properties of the null lines such as denseness are independent of the parametrization, we can consider the
flow of any $\mathcal{X}$--vector field. For instance, we may
choose $X=\lambda_2\partial _{x_1}+\lambda_1\partial _{x_2}$, where we
assume that $\lambda_1,\:\lambda_2>0$. Let the flow of $X$ be
given by $(x_1(t),x_2(t))$. We establish the assertion by computing
the rotation number of this flow (see \cite{ha69} for details). $\lambda_2>0$ implies $\lambda_1dx_1-\lambda_2dx_2=0$.
Since the form $\lambda=\lambda_1dx_1-\lambda_2dx_2$ is closed, we
yield an exact  ordinary differential equation on $\R^2$. Hence we
have to find a $F:\R^2\to\R$ such that $\partial_{x_1}F=\lambda_1$
and $\partial_{x_2}F=-\lambda_2$. The initial condition $x_2(0)$ is determined by $F(0,x_2(0))=c$ for a constant $c$.

Integration of $\partial_{x_1}F=\lambda_1$ yields
$F(x_{1},x_{2})=\int\limits_{0}^{x_{1}}\lambda
_{1}(s,x_{2})ds+f(x_{2})$, where $f^{\prime }(x_2)=-\lambda
_2(x_1,x_2)-\int_0^{x_1}(\partial_2\lambda _1)(s,x_2)ds=-\lambda
_2(x_1,x_2)+\int_0^{x_{1}}\left(\partial _{1}\lambda _{2}\right)
(s,x_{2})ds=-\lambda _{2}(0,x_{2})$. A possible $F$ is given by $F(x_1,x_2)=\int\limits_{0}^{x_{1}}\lambda _{1}(s,x_{2})ds-\int\limits_{0}^{x_{2}}\lambda_{2}(0,s)ds$.
We choose $c=0$ and use the Fourier series of $\lambda_1$ and $\lambda_2$ to get the following equation:
$$ l_1x_1+\sum\limits_{k\not=0 \atop l}\lambda_{1_{kl}}\frac{e^{2\pi ikx_1}-1}{2\pi ik}\cdot e^{2\pi ilx_2}-l_2x_2-\sum\limits_{k \atop l\not=0}\lambda_{2_{kl}}\frac{e^{2\pi ilx_2}-1}{2\pi il}=0. $$
Evaluating in $x_{1}=n\in \Z$ yields $l_1n-l_2x_2(n)-\sum\limits_{k \atop l\not=0}\lambda_{2_{kl}}\frac{e^{2\pi ilx_2}-1}{2\pi il}=0$, hence
$$ \frac{x_2(n)}{n}=\frac{l_1}{l_2}-\frac{1}{n}\cdot\frac{1}{l_2}\sum\limits_{k \atop l\not=0}\lambda_{2_{kl}}\frac{e^{2\pi ilx_2}-1}{2\pi il}=0 $$
Since the Fourier series of smooth functions are absolutely
convergent, we have $\rho=\lim_{n\to
  +\infty}\frac{x_2(n)}{n}=\frac{l_1}{l_2}$, whence the assertion.
\end{proof}
\end{lemma}

We finally get:

\begin{prop}
Let $g$ be conformally equivalent to a closed diagonal metric. Then the following holds:
$1$ and $+\infty$ are the only possible values for $\delta_+^0$. Furthermore, these dimensions are characterized as follows:
\begin{enumerate}
\item $\delta_+^0=1$ iff $\frac{l_1}{l_2}\not\in\Q$ iff all $\mathcal{X}$--null lines are dense.
\item $\delta_+^0=+\infty$ iff $\frac{l_1}{l_2}\in\Q$ iff all $\mathcal{X}$--null lines are closed or asymptotic of a closed $\mathcal{X}$-null line.
\end{enumerate}
\end{prop}

\subsection{$\mu$-surfaces}\label{musurf}

By \ref{const} a half spinor which is harmonic or twistor may be seen
as an object that is constant along the lifts of the corresponding
null lines. Unfortunately, we have no a priori control over the
parallel transport in $Q$, and due to the non--compacteness of
$\tilde{G}(1,1)$, the lifts of the null lines may even be unbounded
(see \ref{denselift}). Therefore, we focus on cases where a direct
link between the null lines on $(M^{1+1},g)$ and the harmonic and twistor spinors can be established, without lifting the null lines to $Q$.

As we did earlier, we restrict our investigation to the case of
positive harmonic spinors. One can prove analogous results by
interchanging suitably $S^+/S^-$ and $\mathcal{X}/\mathcal{Y}$, as
pointed out above.

\begin{defnl}{mu-surface}
Let $\varphi\in \Gamma(S^+)$ be a positive harmonic spinor. A smooth function $\mu_\varphi:M^{1+1}\to \C$ verifying
\begin{enumerate}
\item $\mu_\varphi(x)=0$ if and only if $\varphi(x)=0$ and
\item $\mu_\varphi$ is constant along $\mathcal{X}$-curves
\end{enumerate}
is said to be a {\em mass functional} for $\varphi$. A $\mu$-{\rm surface} is a Lorentzian surface admitting a mass functional for every $\varphi\in \Gamma(S^+)$.
\end{defnl}

We shall give examples of $\mu$--surfaces in the following section
(see \ref{exmusurf1} and \ref{exmusurf2} in conjunction with \ref{muflaeche}). The reason for looking for such mass functionals is the following property:

\begin{lemmal}{heritage}
Let $(M^{1+1},g)$ be a $\mu$--surface. Let $\varphi$ be a positive harmonic spinor and $x\in M$ such that $\varphi(x)=0$. If $(x_n)\subset l$ for a fixed $\mathcal{X}$-line  $l$ converges to $x$, then $\varphi_{|l}\equiv 0$.
\end{lemmal}

The general idea to produce obstructions to the inequality $\delta_+\ge 2$ is to assure that a single zero of a harmonic spinor is propagated along all null lines, therefore forcing the spinor to be zero everywhere. We note the following ''heritage principle'':
\smallskip

\begin{coro}
Let $(M^{1+1},g)$ be a $\mu$--surface. Let $\varphi$ be a positive harmonic spinor, and $l_\infty$ a closed $\mathcal{X}$-line. If $\varphi_{|l_\infty}\equiv 0$, then $\varphi_{|l}\equiv 0$ for every asymptotic $l$ of $l_\infty$.
\end{coro}

Although all Spin structures on a $\mu$--surface can be treated
simultaneously as we shall see, the following proposition illustrates how the
non--trivial Spin structures differ from the trivial one in terms of the parallel transport:

\begin{propl}{denselift}
Let $(M^{1+1},g)$ be a $\mu$--surface with a dense $\mathcal{X}$--line $l$. Then there exists a local section $\tilde{s}:U\to Q$ and a convergent sequence $(x_n)\subset U\cap l$ with $x_n\to x\in U$ such that the following property holds:
If $(a_n)\subset \R$ is defined by $\mathcal{P}^Q_{l:x_0\to x_n}\tilde{s}(x_0)=\tilde{s}(x_n)g_{a_n}$, then for a subsequence $(a_{n_l})$ we have  $a_{n_l}\to \pm\infty$ or
$a_{n_l}\to 0$, that is $\{g_{a_n}\}$ is unbounded in
$\tilde{G}(1,1)$.

\noindent (Recall that according to \ref{spinrep}, $g_{a_n}=\left(\begin{array}{cc}
                                                        a_n & 0 \\
                                                        0   &
                                                        \frac{1}{a_n}
                                                        \end{array}\right)$.)

\begin{proof}
Assume the opposite. Then consider the horizontal lift $l^*$ of $l$ to
$Q$. Extend $l^*$ to a (continuous) section $\tilde{s}_l:M^{1+1}\to Q$ by $\tilde{s}_l(x)=\lim_nl^*(x_n)$ for $x_n\in l\to x$. This limit exists indeed, since $l^*(x_n)=\mathcal{P}^Q_{l:x_0\to x_n}l^*(x_0)$ is bounded in $Q$ by assumption. Hence $Q$ would be isomorphic to the trivial Spin structure.
\end{proof}
\end{propl}

\begin{coro}
If there exists a dense $\mathcal{X}$--line on a $\mu$--surface $(M^{1+1},g)$, then $\delta_+=0$ for every non--trivial Spin structure.

\begin{proof}
Using the notation of the preceding proposition, we have
$$ \varphi(x_n)=[\mathcal{P}^Q_{l:x_0\to x_n} \tilde{s}(x_0),f_{\tilde{s}}(x_0)u_1]=[\tilde{s}(x_n),g_{a_n}f_{\tilde{s}}(x_0)u_1]=[\tilde{s}(x_n),f_{\tilde{s}}(x_0)a_nu_1]\to \varphi(x). $$
Hence, if $a_n\to \pm\infty$, then $f_{\tilde{s}}(x_0)=0$. If $a_n\to 0$, then $\varphi(x)=0$. In both cases, the spinor $\varphi$ has a zero, implying $\varphi\equiv 0$ by \ref{heritage}.
\end{proof}
\end{coro}

From now on, we will mostly consider compact $\mu$-surfaces, though the
techniques and results can be applied to Lorentzian cylinders as
well. Due to the ``denseness obstruction'' \ref{denseob}, we can
restrict our attention to the case where no dense null lines
occur. First, we introduce the subsequent notation:

Let $x\in M^{1+1}$ and $l_1$ and $l_2$ be two closed
$\mathcal{X}$-lines which
do not contain $x$. Since $M^{1+1}$ is homeomorphic to a torus, the
connected components of $M^{1+1}\backslash (l_1\cup l_2)$ are open in
$M^{1+1}$ and homeomorphic to a cylinder without boundary. Let
$C_{l_1l_2}(x)$ denote the cylinder which contains $x$. Its closure is
given by $\overline{C_{l_1l_2}(x)}=C_{l_1l_2}(x)\cup l_1\cup l_2$. In
the case where $l_1=l_2$ as a set, we have
$\overline{C_{l_1l_2}(x)}=M^{1+1}$, so the whole torus itself may be
considered as a closed cylinder. Now let be $x$ such that $l_x$ is an
asymptotic of the two closed null lines $l_1$ and $l_2$. Then the
cylinder $C_{l_1l_2}(x)$ will be written $A_{l_1l_2}(x)$. For further
reference, such a cylinder will be called {\it asymptotic}. Closed
null lines are not allowed to be homotopic to a single point, hence there
are no more closed null lines in any asymptotic cylinder. Since every
asymptotic tends to $l_1$ or $l_2$, we get the

\begin{lemmal}{ascyl}
Let $\varphi$ be a positive harmonic spinor on a compact $\mu$--surface. Then its mass functional $\mu_\varphi$ is constant on every closed asymptotic cylinder.
\end{lemmal}

Thus, if the spinor has a zero in an asymptotic cylinder, it must be zero on the whole cylinder. In order to treat the case where the union of closed null lines is dense in $M^{1+1}$, we introduce a further type of cylinders which does not contain ''ribbons'' of closed null lines:

\begin{defn}
A cylinder $C_{l_1l_2}$ is called {\rm non--resonant} if for any two arbitrary closed null lines $\tilde{l}_1,\tilde{l}_2$ in the closure of $C_{l_1l_2}$ there is an asymptotic $l$ in  $C_{\tilde{l}_1\tilde{l}_2}$.
\end{defn}

\begin{lemmal}{constcyl}
Let $\varphi$ be a positive harmonic spinor, and $C=C_{l_1l_2}$ be a non--resonant cylinder. Then $\mu_\varphi$ is constant on $C$.

\begin{proof}
Consider the set $A:=\{x\in C\;|\;l_x \mbox{ is an asymptotic}\}$. $A$ is open and dense in $C$. As the total differential of $\mu_\varphi$ vanishes on $A$ as a consequence of \ref{ascyl}, the denseness of $A$ implies the result.
\end{proof}
\end{lemmal}

\begin{coro}
Let $(M^{1+1},g)$ be a non--resonant $\mu$--cylinder and let
$\varphi_1,\varphi_2$ be two positive harmonic spinors. If there
exists a $x\in M^{1+1}$ such that  $\varphi_1(x)=\varphi_2(x)$, then
$\varphi_1\equiv\varphi_2$. In particular, every positive harmonic
spinor with a zero is identically zero, and $\delta_+\le 1$.
\end{coro}

\begin{defn}
Let $(M^{1+1},g)$ be a Lorentzian surface and $l$ a closed
$\mathcal{X}$-line. The Spin bundle $Q$ is called $\mathcal{X}${\em
  -trivial along} $l$, if the relation $\mathcal{P}^Q_{l:x\to x}q=q$
holds for every $x\in l$ and $q\in (Q_{|l})_x$.
\end{defn}

\begin{lemma}
Let $\varphi$ be a positive harmonic spinor that has no zero along a closed $\mathcal{X}$-line $l$. Then $Q$ must be $\mathcal{X}$-trivial along $l$.

\begin{proof}
By \ref{parallel} we know that
$\varphi(x)=[q,v]=[\mathcal{P}^Q_{l:x\to x}q,v]$. Hence, if
$\mathcal{P}^Q_{l:x\to x}q=qg$ for a uniquely determined $g\in
Spin(1,1)$, we have $g^{-1}v=v$. It follows $g=id$ by \ref{spinrep}.
\end{proof}
\end{lemma}

So far, we have obtained obstructions to the inequality $\delta_+\ge 2$. In the case where this inequality holds, a third type of cylinder becomes interesting:

\begin{defn}
A closed cylinder $R_{l_1l_2}=\overline{C_{l_1l_2}}$ which does not consist of a single $\mathcal{X}$--line is said to be {\em resonant} if $l_x$ is closed for every $x\in R_{l_1l_2}$.
\end{defn}

\begin{prop}
Let $(M^{1+1},g)$ be a compact $\mu$-surface.
\begin{enumerate}
\item A non--trivial positive harmonic spinor cannot be zero on every resonant cylinder.
\item If $\delta_+\ge 2$, then there exists a $\mathcal{X}$--trivial resonant cylinder $R$ on $M^{1+1}$, that is $Q$ is $\mathcal{X}$-trivial along every closed $\mathcal{X}$--line in $R$.
\end{enumerate}
\begin{proof}
(ii) is a consequence of (i). To prove the first assertion, let us
assume the opposite. It suffices to show that $\mu_{\varphi}$ is
locally constant.

Let $x\in M^{1+1}$. If $l_x$ is an asymptotic, then
$\mu_{\varphi|A_{l_1l_2}(x)}$ is constant by \ref{ascyl}. Otherwise,
$l_x$ is closed. If for every neighbourhood $U$ of $x$ there exists
$x'\in U$ such that $l_x$ is an asymptotic, then $x$ is in the closure
of a non--resonant cylinder $C$. If $x\in int(C)$, then $\mu_\varphi$
is constant on a neighbourhood of $x$ by \ref{constcyl}. If not, then
$x\in\partial C\cap\partial R$, where $R$ is a resonant cylinder. Thus
$\mu_\varphi\equiv 0$ on a neighbourhood of $x$, since
$\mu_{\varphi|C}\equiv const$ and $\mu_{\varphi|R}\equiv 0$ by assumption.
\end{proof}
\end{prop}

On the other hand, whenever there exists an $\mathcal{X}$--trivial resonant
cylinder on $M^{1+1}$, then we can produce harmonic spinors as in
\ref{1con}, since the $\mathcal{X}$--triviality guarantees that the
spinors constructed in this way are well defined. Hence we arrive at the following proposition, generalizing the left--invariant case:

\begin{theoreml}{class}
Let $(M^{1+1},g)$ be a compact $\mu$--surface. Then the only possible dimensions are $\delta_+=0,1$ and $+\infty$. These cases are characterized as follows:
\begin{enumerate}
\item $\delta_+\le 1$ if and only if either \begin{itemize}
        \item there exists a dense $\mathcal{X}$--line in which case we have $\delta_+=0$ for the non--trivial Spin structures or
        \item $M^{1+1}$ is non--resonant or
        \item there exists no $\mathcal{X}$--trivial resonant cylinder on $M^{1+1}$.
        \end{itemize}
\item $\delta_+=+\infty$ if and only if there exists a resonant $\mathcal{X}$--trivial cylinder on $M^{1+1}$. In this case, we have $\delta_+=+\infty$ for every Spin structure.
\end{enumerate}
\end{theoreml}

As we have already seen for the left--invariant case, the dimensions
$\delta_+=0$ and $1$ can occur and depend on the given Spin structure.

\subsection{Spinors and conformal flatness}\label{spinorsandconf}

We will now study the relationship between the existence of harmonic
and twistor spinors and conformal flatness. In particular, we will consider the geometric implications of $\mathcal{X}$--triviality.

As we saw in \ref{confprop}, conformal flatness is related to the
existence of nowhere vanishing time- resp. spacelike conformal vector
fields. With every spinor $\varphi\in\Gamma(S)$ we can canonically associate a vector field that is conformal in the case of a twistor spinor:

\begin{defn}
Let $(M^{p+q},g)$ be an orientable and time--orientable pseudo--Riemannian Spin manifold, and let $\psi\in\Gamma(S)$. We define the {\em associated vector field} $V_\psi$ by the equation
\begin{equation*}
g\left( V_{\psi },W\right) =i^{p+1}\left\langle W\cdot \psi ,\psi
\right\rangle
\end{equation*}
for $W\in \mathfrak{X}\left( M^{p+q}\right) $.
\end{defn}

A direct computation yields the following proposition (see, for
instance, \cite{ba99}):

\begin{propl}{conform}
Let $\varphi \in \Gamma\left( S\right)$ be a twistor spinor. Then $V_{\psi }$ is a conformal vector field. More precisely, we have $\mathcal{L}_{V_{\psi }}g=%
\frac{4}{n}Re\left( i^{p+1}\left\langle D\psi ,\psi
\right\rangle \right) g$.
\end{propl}

Next, we determine the associated vector field of a spinor in signature $(1,1)$:

\begin{lemmal}{spinvf}
Let $\left( M^{1+1},g\right)$ be time--orientable and $\psi \in
\Gamma\left( S\right)$. Let $s=\left(s_1,s_2\right):U\to P$ be an
orthonormal frame with a lift $\tilde{s}$ to $\tilde{Q}_{\xi}$
(cf. section \ref{prspin}). Let
$\psi_{\tilde{s}}=\psi^+_{\tilde{s}}u_1+\psi^-_{\tilde{s}}u_{-1}\in
C^\infty\left(U,\Delta_{1,1}\right)$ be the local trivialization of
$\psi$ with respect to $\tilde{s}$. Then
$V_\psi=\left|\psi^+_{\tilde{s}}\right|^2X_s-\left|\psi^-_{\tilde{s}}\right|^2Y_s$.
In particular, $V_\psi$ is a causal vector feld which is timelike if the local components $\psi^+_s$ and $\psi^-_s$ have no zeros, and lightlike in case of a half--spinor without zeros.

\begin{proof}
Let  $w_1$ and $w_2$ be the local components of $W\in
\mathfrak{X}\left( M\right) $ with respect to $s$, that is
$W=w_{1}s_{1}+w_{2}s_{2}=\left[ s,w_{1}e_{1}+w_{2}e_{2}\right] $. A
direct computation yields $<W\cdot \psi ,\psi
>=(\left|\psi^-_{\tilde{s}}\right| ^{2}+\left|\psi^+_{\tilde{s}}
\right|^{2})w_{1}+(\left|\psi^+_{\tilde{s}}\right| ^{2}-\left|
  \psi^-_{\tilde{s}}\right| ^{2})w_{2}$. If $V_{\psi }=V_{\psi
  1}s_1+V_{\psi 2}s_2$, we get $V_{\psi 1}=\left
  (\left|\psi^+_{\tilde{s}} \right| ^{2}+\left|
    \psi^-_{\tilde{s}}\right| ^{2}\right )$ and $V_{\psi 2}=\left (\left| \psi^+_{\tilde{s}}\right| ^{2}-\left|
\psi^-_{\tilde{s}}\right| ^{2}\right )$. Furthermore, since $\lambda=g(X_s,Y_s)>0$ and $g\left( V_{\psi },V_{\psi }\right) =-2\lambda \left| \psi^+_{\tilde{s}}\right| ^{2}\left| \psi^-_{\tilde{s}}\right| ^{2}$, the vector field $V_{\psi }$ is causal.
\end{proof}
\end{lemmal}

\begin{rem}
Let $X$ be a $\mathcal{X}$--vector field and $\varphi$ a positive
harmonic spinor. Then $X\cdot\varphi=0$. In particular we get $V_{\varphi} \cdot \varphi=0$.
\end{rem}

\begin{corol}{confflattw}
Let $\left( M^{1+1},g\right)$ and $\psi ^\pm\in\Gamma(S^\pm)$ be two twistor spinors without zeros. Then $(M^{1+1},g)$ is conformally flat. In particular if $M$ is compact, then $(M^{1+1},g)$ is complete.
\end{corol}

\begin{rem}
According to \ref{correspondence} the same result holds for harmonic instead of twistor spinors.
\end{rem}

As we have already seen for $\mathcal{X}$-triviality, harmonic
resp. twistor spinors without zeros induce certain ``flatness''
properties. Therefore, we will look more closely to Lorentzian
surfaces that admit nowhere vanishing solutions to the harmonic
resp. twistor equation.

First, we recall the following statement:

\begin{propl}{divharm}
Let $\varphi\in\Gamma(S)$ be a harmonic spinor on $(M^{p+q},g)$. Then $div(V_\varphi)=0$.
\end{propl}

\begin{defnl}{scfdef}
A time--orientable Lorentzian surface $(M^{1+1},g)$ is said to be $\mathcal{X}-$ resp. $\mathcal{Y}${\em -conformally flat} if there exists a global orthonormal frame $s=(s_1,s_2)$ such that $div(X_s)=0$ resp. $div(Y_s)=0$. We call a Lorentzian surface $(M^{1+1},g)$ {\em semi-conformally flat (s.c.f.)} if $(M^{1+1},g)$ is either $\mathcal{X}$-- or $\mathcal{Y}$--conformally flat.
\end{defnl}

The notion of semi-conformal flatness will be justified in \ref{xyconfflat}.

As we did for spinors, we will concentrate on $\mathcal{X}$--conformally flat surfaces; analogous statements hold for $\mathcal{Y}$--conformally flat ones.

Since every $\mathcal{X}$--vector field can be written as $X_s$ with
respect to a suitably chosen basis, a time--orientable Lorentzian
surface is $\mathcal{X}$--conformally flat if and only if there exists
a $\mathcal{X}$--vector field $X$ such that $div(X)=0$. Furthermore,
it follows that the notion of semi-conformal flatness is invariant
under conformal change of the metric: If there exists a
$\mathcal{X}$-vector field $X$ on $(M^{1+1},g)$ with $div(X)=0$, then
we can find another $\mathcal{X}$-vector field
$\tilde{X}\in\mathfrak{X}(M)$ with $\widetilde{div}(\tilde{X})=0$
(where $\widetilde{div}$ denotes the divergence operator associated
with the conformally changed metric $\tilde{g}=\lambda g$) as can be
seen from the formula $\widetilde{div}(V)=V(\ln(\lambda))+div(V)$.

\begin{exmsl}{exmusurf1}
(i) As we remarked in section \ref{cldiagmetrics}, closed diagonal
metrics admit lightlike divergence-free vector fields and are
therefore s.c.f..

\noindent (ii) We are going to exhibit further examples by a direct
computation of the divergence: Let us consider a Lorentzian torus with
standard coordinates $(x_1,x_2)$ and volume form $\omega$. Let
$V=k\partial_{x_1}+l\partial_{x_2}$. Using the formula
$d(i_X\omega)=div(X)\omega$, we get $div(V)= \partial_1 k+ \partial_2
l +\frac{1}{2}V(\ln|\det(g)|)$. In particular, if $\det(g)\equiv 1$,
then $div(V)= \partial_1 k+ \partial_2 l$. As semi-conformal flatness is a conformal invariant, we may always assume -- by rescaling the metric with the factor $-\frac{1}{det(g)}$ -- this assumption to be fulfilled.
For instance, if we consider the family of metrics given by
(\ref{san97g1}) and (\ref{san97g2}), we get:
\end{exmsl}

\begin{corol}{exmusurf2}
Every metric in $\mathcal{G}_1\cup\mathcal{G}_2$ defines a s.c.f. surface.
\end{corol}

\begin{propl}{hcfchar}
Let $(M^{1+1},g)$ be a time--orientable Lorentzian surface. Then the following assertions are equivalent:
\begin{enumerate}
\item $(M^{1+1},g)$ is $\mathcal{X}$--conformally flat
\item There exists a global section $M^{1+1}\to P$ such that $div(X_s)=0$
\item There exists a $\mathcal{X}$--vector field $X$ such that $div(X)=0$
\item $\nabla_{X_s}^{LC}s_i=0$ for $i=1,2$
\item There exists a positive harmonic spinor without zeros
\item There exists a negative twistor spinor without zeros
\end{enumerate}
\begin{proof}
Only the implications (ii) $\Rightarrow$ (iv), (iv) $\Rightarrow$ (ii), (v) $\Rightarrow$ (iii) and (ii)
$\Rightarrow$ (v) need proof.

\noindent (ii) $\Leftrightarrow$ (v): By a direct application of the
Koszul formula, we prove the

\begin{lemmal}{divequiv}
$$ -\frac{1}{2}g(X_s,[X_s,Y_s])=div(X_s)=g(\nabla^{LC}_{X_s}s_1,s_2)=-g(\nabla^{LC}_{X_s}s_2,s_1) $$
\end{lemmal}

Then (\ref{spinder11}) implies:

\begin{corol}{spinder}
Locally, we have the identity
$\nabla^S_{X_s}[\tilde{s},\varphi_{\tilde{s}}]=[\tilde{s},X_s(\varphi_{\tilde{s}})-\frac{1}{2}div(X_s)e_1\cdot
e_2\cdot \varphi_{\tilde{s}}]$. In particular, if $(M^{1+1},g)$ is $\mathcal{X}$--conformally flat, we get $\nabla^S_{X_s}[\tilde{s},\varphi_{\tilde{s}}]=[\tilde{s},X_s(\varphi_{\tilde{s}})]$.
\end{corol}

Furthermore, $g(\nabla^{LC}_{X_s}s_i,s_i)=0$ for $i,j\in\{1,2\}$. By the lemma, we have $g(\nabla^{LC}_{X_s}s_1,s_2)=-g(\nabla^{LC}_{X_s}s_2,s_1)=div(X_s)$, whence the equivalence.

(v) $\Rightarrow$ (iii):  Since $div(e^fX)=e^f(X(f)+div(X))$, the
$\mathcal{X}$--vector field $e^fX$ will be divergence--free if and
only if $X(f)=-div(X)$ holds. By \ref{spinvf}, we have
$V_\varphi=\lambda X_s$ with $\lambda\not= 0$. Application of \ref{divharm} yields $X_s(\ln|\lambda|)=-div(X_s)$.

(ii) $\Rightarrow$ (v):  Let $\tilde{s}:M^{1+1}\to Q_0$ be a global
section in the trivial bundle. Then $\varphi_{\tilde{s}}=u_1$ defines a positive harmonic spinor without zeros.
\end{proof}
\end{propl}

\begin{paral}{deltage1}{\sc Corollary 1.} {\it
On a $\mathcal{X}$--conformally flat surface, we have $\delta_+^0\ge 1$.
}
\end{paral}

\begin{paral}{xyconfflat}{\sc Corollary 2.} {\it
$(M^{1+1},g)$ is conformally flat if and only if $(M^{1+1},g)$ is $\mathcal{X}$-- and $\mathcal{Y}$--conformally flat.
}

\begin{proof}
For a flat metric, every constant defines a harmonic resp. twistor
spinor with respect to the trivial Spin structure. The implication
follows then from \ref{confspininv}. We yield the converse from \ref{confflattw}.
\end{proof}
\end{paral}

Next, we will prove some further properties of $\mathcal{X}$--conformally flat surfaces:

\begin{defn}
A Lorentzian surface $(M^{1+1},g)$ is said to be $\mathcal{X}$-- resp. $\mathcal{Y}$--{\em complete}, if every $\mathcal{X}$-- resp. $\mathcal{Y}$--geodesic is complete. A Lorentzian surface that is $\mathcal{X}$-- resp. $\mathcal{Y}$--complete is said to be {\em semi--null complete}.
\end{defn}

\begin{prop}
A compact $\mathcal{X}$--conformally flat surface is $\mathcal{X}$--complete.

\begin{proof}
By \ref{hcfchar} (v), we have $\nabla^{LC}_{X_s}X_s=0$, so that the $X_s$--geodesics are given by the flow of $X_s$.
\end{proof}
\end{prop}

\begin{exms}

\noindent (i) The interdependence of semi-conformal flatness, positive harmonic spinors without zeros and semi-completeness is demonstrated by \ref{analex}: Since the harmonic half--spinors we found have no zeros, the surface must be semi-complete. But we showed that this metric is in $\mathcal{G}_2$, so it is not conformally flat. Hence, from \ref{xyconfflat} follows that $(T^{1+1},g)$ is complete for one type of isotropic geodesics, and that there must be incomplete geodesics for the other type, in accordance with proposition \ref{nullcomp}.

\noindent (ii) Consider the following example taken from
\cite{rosa93}: Let $\tau:[0,1]\to \R$ be a smooth function with
$\tau(a)=0$, but $\tau'(a)\not=0$, and whose support is strictly
contained in $[0,1]$. Extend $\tau$ periodicly on the whole real line
and define $g^{\tau}_{(x,y)}=2dxdy-\tau(x)dy^2$. Then
$g^{\tau}\in\mathcal{G}'$, hence $g^{\tau}$ is lightlike
incomplete. For instance,
$\gamma(t)=(a,\tau'(a)\ln(t+\frac{1}{\tau'(a)}))$ is a closed
incomplete geodesic which  without loss of generality we assume to be
$\mathcal{X}$. Thus any positive harmonic spinor must be zero on
$\gamma$ . But as $\tau_{|[-\epsilon,\epsilon]}\equiv 0$ for
$\epsilon$ sufficiently small, $(T^{1+1},g)$ contains an
$\mathcal{X}$-trivial resonant cylinder, and therefore
$\delta_+=+\infty$. This example shows that there exists not conformally
flat tori with $\mathcal{X}$-trivial resonant cylinder which are not
$\mathcal{X}$-conformally flat.
\end{exms}

\begin{prop}
On a $\mathcal{X}$--conformally flat surface, we have $\mathfrak{H}_+\cong \mathfrak{T}_-$.

\begin{proof}
The maps $\Phi:\mathfrak{H}_+\to \mathfrak{T}_-,\; \varphi\mapsto
\frac{i}{2}Y_s\cdot \varphi$ and $\Psi:\mathfrak{T}_-\to
\mathfrak{H}_+,\; \psi\mapsto \frac{i}{2}X_s\cdot \varphi$ are bundle
isomorphisms inverse of one another.
\end{proof}
\end{prop}

\begin{propl}{muflaeche}
Every s.c.f. surface is a $\mu$--surface.

\begin{proof}
Let $s=(s_1,s_2):M\to P$ be a (global) orthonormal frame with $div(X_s)=0$. We define
$$\mu_\varphi(x)=\langle Y_s\cdot\varphi,\varphi\rangle. $$
Let $\tilde{s}$ be a local lift of $s$ to
$\tilde{Q}_{s_1}$. For this section, let $\varphi=[\tilde{s},\varphi_{\tilde{s}}]$. Then
$\mu_\varphi(x)=\langle Y_s\cdot\varphi,\varphi\rangle(x)=-2|\varphi_{\tilde{s}}(x)|^2$, where
$|\cdot|$ denotes the absolute value function on $\C$. Thus (i) and
(ii) of definition \ref{mu-surface} hold.

Since $\varphi$ is a positive harmonic spinor, we have
$\nabla^S_{X_s}\varphi=0$. Consequently, we get
$X_s(\mu_\varphi)=\langle\nabla^{LC}_{X_s} Y_s\cdot\varphi,\varphi \rangle.$
Since $\nabla^{LC}_{X_s} Y_s=-\nabla^{LC}_{X_s} s_1+\nabla^{LC}_{X_s} s_2=0$ by \ref{hcfchar} (iv), the assertion follows.
\end{proof}
\end{propl}

The classification of the possible dimensions of $\delta_+$ in \ref{class} may be restated as follows:

\begin{theoreml}{scfthm}
Let $(M^{1+1},g)$ be a compact $\mathcal{X}$--conformally flat
Lorentzian surface. Then $\delta_+=\tau_-$ and the only possible dimensions for $\delta_+$ are $0,1$ and $+\infty$. These cases are characterized as follows:
\begin{enumerate}
\item $\delta_+\le 1$ if and only if either \begin{itemize}
        \item there exists a dense $\mathcal{X}$--line in which case we have $\delta_+=0$ for the non--trivial Spin structures, or
        \item $M^{1+1}$ is non-resonant or
        \item there exists no $\mathcal{X}$--trivial resonant cylinder on $M^{1+1}$.
        \end{itemize}
Furthermore, $\delta_+^0=1$.
\item $\delta_+=+\infty$ if and only if there exists an $\mathcal{X}$--trivial resonant cylinder on $M^{1+1}$. In this case we have $\delta_+=+\infty$ for every Spin structure.
\end{enumerate}
\end{theoreml}

Next we will show that in some sense $\mathcal{X}-$conformal flatness
is forced by positive harmonic spinors which have non--zero ''mass'':

\begin{defn}
Let $(M^{1+1},g)$ be a Lorentzian surface and $(P,\pi,M^{1+1};G)$ a
principal fiber bundle over $M^{1+1}$. A (local) section $s:U\to P$ is
said to be $\mathcal{X}$-- resp. $\mathcal{Y}$--{\em parallel} if for
every $\mathcal{X}$-- resp. $\mathcal{Y}$--curve $\alpha:[a,b]\to U$, we have $\mathcal{P}^P_{\alpha:a\to b} s(\alpha(a))=s(\alpha(b))$, where $\mathcal{P}_\alpha^P$ denotes the parallel transport in $P$ along $\alpha$.
\end{defn}

\begin{lemmal}{divpara}
Let $s=(s_1,s_2):U\to P$ be a local section of the orthonormal frame bundle $P$. Then $s$ is $\mathcal{X}$-parallel if and only if $div(X_s)=0$. Furthermore, if $s$ can be lifted to a section $\tilde{s}:U\to Q$ of $Q$, then $\tilde{s}$ is $\mathcal{X}$-parallel if and only if $div(X_s)=0$.

\begin{proof}
Let $\alpha$ be the flow generated by $X_s$ in $U$, and
let $\mathcal{P}^{LC}_\alpha$ be the usual parallel transport in
$TM^{1+1}$ along $\alpha$ induced by $\mathcal{P}^P$. We have
$\nabla^{LC}_{X_s}s_j(x)=\frac{d}{dt}\mathcal{P}^{LC}_{\alpha:t\to
  0}s_j(\alpha(t))_{|_{t=0}}=\frac{d}{dt}[s(\alpha(0)),e_j]=0$, hence $div(X_s)=0$ by \ref{hcfchar} (v).

For the converse, let $Z$ denote the Levi-Civita-connection in $P$. Since
$div(X_s)=-\frac{1}{2}g(\nabla^{LC}_{X_s}s_1,s_2)=0$, we get
$s^*Z(X)=Z(ds(X))=0$ for every $\mathcal{X}$-vector field $X$ (cf. \ref{spinder}). Hence $s^*Z(\alpha'(t))=0$ for every $\mathcal{X}$-curve $\alpha:[a,b]\to M^{1+1}$, that is $\alpha^*_{s(\alpha(a))}=$ lift of $\alpha$ starting in $s(\alpha(a))=s\circ\alpha$. It follows that $\mathcal{P}^P_\alpha s=\alpha^*_{s(\alpha(a))}=s(\alpha(b))$.

Since $f\circ\mathcal{P}^Q=\mathcal{P}^P\circ f$ and
$\tilde{s}^*\tilde{Z}(X_s)=-div(X_s)\omega =0$ for the lifts of $s$
and $Z$ to $Q$, we deduce the same result for the Spin bundle $Q$.
\end{proof}
\end{lemmal}

The notion of $\mathcal{X}$--triviality can then be reformulated as follows:

\begin{coro}
Let $l$ be a closed $\mathcal{X}$-line. Then $Q$ is $\mathcal{X}$-trivial along $l$ if and only if $l$ can be parametrized such that $div(l')=0$. In particular, such a parametrization makes $l$ into a geodesic.
\begin{proof}
Choose an orthonormal frame $s$ such that $l'=s_1+s_2$ and repeat the reasoning of \ref{divpara}.
\end{proof}
\end{coro}

\begin{rem}
As in the case of conformal flatness the condition $div(X)=0$ can
always be locally realized: Indeed, if $\beta:(a,b)\to U$ is a
$\mathcal{Y}$-curve, pick a section $s:|\beta|\to P$ and extend this
section on $U$ by parallel transport of $s(\beta(t))$ in the $\mathcal{X}$-direction. Hence $s:U\to P$ is $\mathcal{X}$-parallel by construction and well defined if $U$ is conveniently chosen, i.e. $\beta$ intersects every $\mathcal{X}$-line in $U$ only once and there are no closed $\mathcal{X}$-lines (e.g. if $U$ is simply connected). Then $div(X_s)=0$ by \ref{divpara}.
\end{rem}

\subsection{Conclusion} \label{conclusion}

Like for compact Riemannian surfaces, $\delta_+$ depends both on the
conformal class of the metric and on the Spin structure. $\delta_+$
may be unbounded, in contrast to what is known for the Riemannian
case, where the dimension is bounded by $[\frac{g+1}{2}]$ ($g$
denoting the genus of the surface) -- see \cite{hi74}. Furthermore, for
Lorentzian surfaces we have a certain symmetry between harmonic and twistor spinors.

In the case of $\mu$--surfaces, the conformal invariants $\delta$ and
$\tau$ reflect the global behaviour of the null lines. In some regular
cases, where the global and local behaviour is quite similar (e.g. for
simply connected surfaces or resonant tori), $\delta$ and $\tau$ are
$+\infty$. If a ''pathological'' behaviour such as dense null lines
occurs, then $\delta$ and $\tau$ are less than or equal $2$, and we
have a kind of ''dynamic'' dependence on the conformal class. No
intermediate values are attained. Although $\delta$ and $\tau$ are
weaker conformal invariants than the null lines, in some cases they allow us to
distinguish between conformal classes. Furthermore,
solutions with ''mass'', that is solutions without zeros, force
conformal flatness.  All techniques used -- above all the
characterization of harmonic and twistor half spinors as a kind of parallel spinors along the lightlike distributions -- are
genuine for the signature $(1,1)$. On the other hand the case of a
pseudo-Riemannian signature $(p,q)$ with $p+q\ge 3$ is significantly different. For
instance, the dimension of the space of twistor spinors on a connected
pseudo-Riemannian manifold is bounded by $2^{\left[\frac{p+q}{2}\right]+1}$ (see \cite{ba99}
and \cite{bfgk91})

It is not clear altogether to what extent these techniques can be
applied to a general Lorentzian surface or which are the geometric
obstructions for doing so. The next obvious step would be to
investigate the class of asymptotic cylinders. One could try to find
counterexamples of the ''heritage principle'', that is, harmonic or
twistor spinors which are zero on the closed null lines, but have no
zeros on the asymptotic cylinder itself, or metrics for which $\delta_{\pm}$ or $\tau_{\pm}$ may attain values other than $0,1$ or $+\infty$.


\begin{thebibliography}{99}
        \bibitem[Ba81]{ba81} {\sc H.Baum.} {\it Spin-Strukturen und Dirac--Operatoren \"{u}ber pseudoriemannschen Mannigfaltigkeiten.} {\rm Teubner, Leipzig, 1981.}
        \bibitem[Ba99]{ba99} {\sc H.Baum.} {\it Lorentzian twistor
        spinors and CR-geometry.} {\rm Diff.~geom.~Appl. {\bf 11}
        (1999), 69-96.}
        \bibitem[BFGK91]{bfgk91} {\sc H.Baum, T.Friedrich,
        R.Grunewald, I.Kath.} {\it Twistor and Killing Spinors on
        Riemannian manifolds.} {\rm Teubner-Texte zur Mathematik 124,
        Teubner, Stuttgart, 1991.}
        \bibitem[B\"aSc92]{basc92} {\sc C.B\"ar, P.Schmutz.} {\it Harmonic spinors on Riemannian surfaces.} {\rm Ann.~Glob.~Anal.~Geom. {\bf 10} (1992), 263-273.}
        \bibitem[Ha69]{ha69} {\sc J.Hale.} {\it Ordinary differential equations.} Wiley-Interscience, New York, 1969.
        \bibitem[Hi74]{hi74} {\sc N.Hitchin.} {\it Harmonic Spinors.} Adv.~in Mathematics, \textbf{14} (1974), 1-55.
        \bibitem[Ka68]{ka68} {\sc M.Karoubi.} {\it Alg\`{e}bres de
Clifford et K-Th\'{e}orie.} {\rm Ann. scient. \'{E}c.~Norm.~Sup., 4$^{\text{e }}$s\'{e}rie, \textbf{1} (1968), 161-270.}
        \bibitem[RoS\'a93]{rosa93} {\sc A.Romero, M.S\'{a}nchez.} {\it New properties and examples of incomplete Lorentzian tori.} J.~Math.~Phys. {\bf 35 }(1994), 1992-1997.
        \bibitem[S\'a97]{sa97} {\sc M.S\'{a}nchez.} {\it Structure of Lorentzian tori with a Killing vector field.} Trans.~Am.~Math.~Soc. {\bf 349,} No.3 (1997), 1063-1080.
        \bibitem[We96]{we96} {\sc T.Weinstein.} {\it An introduction to Lorentz surfaces.} De Gruyter, Berlin, 1996.
\end{thebibliography}
\end{document}